\documentclass[12pt]{amsart}
\usepackage{amssymb,amsthm,amsmath,a4}
\newtheorem{thm}{Theorem}[section]
\newtheorem{theorem}[thm]{Theorem}
\newtheorem{corollary}[thm]{Corollary}
\newtheorem{lemma}[thm]{Lemma}

\theoremstyle{definition}
\newtheorem{definition}[thm]{Definition}
\newtheorem{example}[thm]{Example}

\theoremstyle{remark}
\newtheorem{remark}[thm]{Remark}
\numberwithin{equation}{section}

\newcommand{\eps}{{\varepsilon}}

\def\CC{{\mathcal C}}
\def\IC{{\mathcal I}}
\def\JC{{\mathcal J}}
\def\KC{{\mathcal K}}
\def\LC{{\mathcal L}}
\def\MC{{\mathcal M}}
\def\NC{{\mathcal N}}
\def\OC{{\mathcal O}}
\def\PC{{\mathcal P}}
\def\XC{{\mathcal X}}
\def\YC{{\mathcal Y}}
\def\VC{{\mathcal V}}

\def\DR{{\mathrm D}}
\def\NR{{\mathrm N}}
\def\Lip{\mathrm{Lip}}
\def\lip{\mathrm{lip}}
\def\Osc{\mathrm{Osc}\,}

\def\MB{\mathbf M}
\def\PB{\mathbf P}
\def\VB{\mathbf V}
\def\kB{\mathbf k}

\def\bR{\mathrm b}
\def\dR{\mathrm d}
\def\eR{\mathrm e}

\newcommand{\R}{\mathbb{R}}
\newcommand{\N}{\mathbb{N}}

\newcommand{\C}{\mathbb{C}}
\newcommand{\Z}{\mathbb{Z}}
\newcommand{\supp}{\operatorname{supp}}
\newcommand{\diam}{\operatorname{diam}}
\newcommand{\codim}{\operatorname{codim}}
\newcommand{\dist}{\operatorname{dist}}

\def\leq {\leqslant}
\def\geq {\geqslant}

\begin{document}

\title[Weil formula]{Weil asymptotic formula
for the Laplacian on domains with rough boundaries}

\author[Netrusov]{Yu. Netrusov\,$^1$}
\author[Safarov]{Yu. Safarov}
\address
{School of Mathematics, University of Bristol, Bristol BS8 1TW,
UK} \email{Y.Netrusov@bris.ac.uk}
\address
{Department of Mathematics, King's College, London WC2R 2LS, UK}
\email{ysafarov@mth.kcl.ac.uk}
\date{October 2003}
 \footnote{Research supported by EPSRC
grant GR/A00249/01}

\begin{abstract}
We study asymptotic distribution of eigenvalues of the Laplacian
on a bounded domain in $\,\R^n$. Our main results include an
explicit remainder estimate in the Weyl formula for the Dirichlet
Laplacian on an arbitrary bounded domain, sufficient conditions
for the validity of the Weyl formula for the Neumann Laplacian on
a domain with continuous boundary in terms of smoothness of the
boundary and a remainder estimate in this formula. In particular,
we show that the Weyl formula holds true for the Neumann Laplacian
on a $\,\Lip_\alpha$-domain whenever $\,(d-1)/\alpha<d\,$, prove
that the remainder in this formula is
$\,O(\lambda^{(d-1)/\alpha})\,$ and give an example where the
remainder estimate $\,O(\lambda^{(d-1)/\alpha})\,$ is order sharp.
We use a new version of variational technique which does not
require the extension theorem.
\end{abstract}

\maketitle

\numberwithin{equation}{section}

\section*{Introduction}

Let $-\Delta_\NR$ be the Neumann Laplacian on a bounded domain
$\,\Omega\subset\R^d\,$ and $N_\NR(\Omega,\lambda)$ be the number
of its eigenvalues which are strictly smaller than $\lambda^2$; if
the number of these eigenvalues is infinite or $-\Delta_\NR$ has
essential spectrum below $\lambda$ then we define
$N_\NR(\Omega,\lambda):=+\infty$. Similarly, let $-\Delta_\DR$ be
the Dirichlet Laplacian on $\Omega$ and $N_\DR(\Omega,\lambda)$ be
the number its eigenvalues lying below $\lambda^2$. We shall omit
the lower index $\,\DR\,$ or $\,\NR\,$ and simply write
$\,\Delta\,$ or $\,N(\Omega,\lambda)\,$ if the corresponding
statement refers both to the Dirichlet and Neumann Laplacian.

According to the Weyl formula,
\begin{equation}\label{0.1}
N(\Omega,\lambda)\ -\ C_{d,W}\,\mu_d(\Omega)\,\lambda^d\ =\
o(\lambda^d)\,,\qquad\lambda\to+\infty\,,
\end{equation}
where $\,\mu_d(\Omega)\,$ is the $\,d$-dimensional Lebesgue
measure of $\,\Omega\,$ and $\,C_{d,W}\,$ is the Weyl constant
(see Subsection \ref{S1.1}). If $\,N=N_\DR\,$ then the Weyl
formula holds for all bounded domains \cite{BS}. If, in addition,
the upper Minkowski dimension of the boundary is equal to
$\,d_1\in(d-1,d)\,$ then
\begin{equation}\label{0.2}
N(\Omega,\lambda)\ -\ C_{d,W}\,\mu_d(\Omega)\,\lambda^d\ =\
O(\lambda^{d_1})\,,\qquad\lambda\to+\infty\,.
\end{equation}
The asymptotic formula (\ref{0.2}) with $\,N=N_\DR\,$ is well
known and is proved in many papers, for instance, in \cite{BLi}
and \cite{Sa} where the authors obtained estimates with explicit
constants. This formula remains valid for the Neumann Laplacian
whenever $\,\Omega\,$ has the extension property (see below). Note
that $\,d_1\,$ may well coincide with $\,d\,$, in which case
(\ref{0.2}) is useless.

If $\,N=N_\NR\,$ then (\ref{0.1}) is true only for domains with
sufficiently regular boundaries. In the general case $\,N_\NR\,$
does not satisfy (\ref{0.1}); moreover, the Neumann Laplacian on a
bounded domain may have a nonempty essential spectrum (see, for
example, \cite{HSS} or \cite{Si}). The necessary and sufficient
conditions for the absence of the essential spectrum in terms of
capacities have been obtained in \cite{M1}. In \cite{BD} the
authors proved that $\,N_\NR(\Omega,\lambda)\,$ is polynomially
bounded whenever the Sobolev space $\,W^{1,2}(\Omega)\,$ is
embedded in $\,L^q(\Omega)\,$ for some $\,q>2\,$. If the
log-Sobolev inequality holds on $\,\Omega\,$ then
$\,N_\NR(\Omega,\lambda)\,$ is exponentially bounded \cite{Ma}.

For domains $\,\Omega\,$ with sufficiently smooth boundaries,
(\ref{0.1}) is true for the both functions $\,N_\DR\,$ and
$\,N_\NR\,$ and the remainder (i.e., the right hand side) is
$\,O(\lambda^{d-1})\,$ \cite{Iv1}, \cite{Se}. The proof is based
on the study of propagation of singularities for the corresponding
evolution equation (see \cite{Iv3} or \cite{SV}). If $\,\Omega\,$
has a rough boundary then the propagation of singularities near
$\,\partial\Omega\,$ cannot be effectively described and one has
to invoke the variational technique.

Let $\,\Omega_\delta^\bR\,$ and $\,\Omega_\delta^\eR\,$ be the
internal and external $\,\delta$-neighbourhoods of
$\,\partial\Omega\,$ respectively. The classical variational proof
of the Weyl formula involves covering the domain by a finite
collection of disjoint cubes $\,\{Q_j\}_{j\in\JC}\,$ and using the
Dirichlet--Neumann bracketing. It is convenient to assume that
$\,\{Q_j\}_{j\in\JC}\,$ is the subset of the family of Whitney
cubes covering $\,\Omega\bigcup\Omega_\delta^\eR\,$ (see
Theorem~\ref{T3.3}), which consists of the cubes $\,Q_j\,$ such
that $\,Q_j\bigcap\Omega\ne\emptyset\,$.

In view of the Rayleigh--Ritz variational formula, we have the
estimates $\,\sum_{j\in\JC_0}N_\DR(Q_j,\lambda)\leq
N_\DR(\Omega,\lambda)\leq \sum_{j\in\JC}N_\NR(Q_j,\lambda)\,$,
where $\,\{Q_j\}_{j\in\JC_0}\,$ is the set of cubes $\,Q_j\,$
lying inside $\,\Omega\,$. If $\,\mu_d(\partial\Omega)=0\,$ then,
estimating $\,N_\DR(Q_j,\lambda)\,$ and $\,N_\NR(Q_j,\lambda)\,$
for each $\,j\,$ and taking $\,\delta=\lambda^{-1}\,$, we obtain
(\ref{0.1}) and (\ref{0.2}) for the Dirichlet Laplacian. It is
possible to get rid of the condition $\,\mu_d(\partial\Omega)=0\,$
but this requires additional arguments.

Similarly, the Rayleigh--Ritz formula implies that
\newline
$\,\sum_{j\in\JC_0}N_\DR(Q_j,\lambda)\leq
N_\NR(\Omega,\lambda)\leq\sum_{j\in\JC_{m\delta}}N_\NR(Q_j,\lambda)+
N_\NR(\bigcup_{j\in\JC\setminus\JC_{m\delta}}Q_j\bigcap\Omega,\lambda)\,$,
where $\,\{Q_j\}_{j\in\JC_{m\delta}}\,$ is the set of cubes lying
inside $\,\Omega\setminus\Omega_{m\delta}^\bR\,$. If for some
$\,m\in\N\,$ and all sufficiently small positive $\,\delta\,$
there exist uniformly bounded extension operators from the Sobolev
space $\,W^{1,2}(\Omega_{m\delta}^\bR)\,$ to
$\,W^{1,2}(\Omega_{m\delta}^\bR\bigcup\Omega_\delta^\eR)\,$ then
$\,N_\NR(\bigcup_{j\in\JC\setminus\JC_{m\delta}}Q_j\bigcap\Omega,\lambda)\leq
N_\NR(\bigcup_{j\in\JC\setminus\JC_{m\delta}}Q_j,C\lambda)
=\sum_{j\in\JC\setminus\JC_{m\delta}}N_\NR(Q_j,C\lambda)\,$, where
$\,C\,$ is a sufficiently large constant. If, in addition,
$\,\mu_d(\partial\Omega)=0\,$ then, estimating the counting
functions on the cubes and taking $\,\delta=\lambda^{-1}\,$, we
obtain (\ref{0.1}) and (\ref{0.2}) for
$\,N_\NR(\Omega,\lambda)\,$.

However, the known extension theorems require certain regularity
conditions on the boundary (for instance, it is sufficient to
assume that $\,\partial\Omega\,$ belongs to the Lipschitz class or
satisfies the cone condition). Domains with very irregular
boundaries do not have the $\,W^{1,2}$-extension property, in
which case the above scheme does not work the Neumann Laplacian.
To the best of our knowledge, in all papers devoted to the Weyl
formula for $\,N_\NR(\Omega,\lambda)\,$ the authors either
implicitly assumed that the domain has the $\,W^{1,2}$-extension
property or directly applied a suitable extension theorem.

The main aim of this paper is to introduce a different technique
which does not use an extension theorem. Instead of disjoint
cubes, we cover the domain $\,\Omega\,$ by a family of relatively
simple sets $\,S_m\subset\Omega\,$. For each of these sets the
counting function $\,N(S_m,\lambda)\,$ can be effectively
estimated from below and above. The sets $\,S_m\,$ may overlap
but, under certain conditions on $\,\Omega\,$, the multiplicity of
their intersection does not exceed a constant depending only on
the dimension $\,d\,$.

This allows us to apply the Dirichlet--Neumann bracketing and
obtain the Weyl asymptotic formula with a remainder estimate for
the Neumann Laplacian on domains without the extension property
(Theorem \ref{T1.3}). The remainder term in this formula may well
be of higher order than the first term. Then our asymptotic
formula turns into an estimate for $\,N_\NR(\Omega,\lambda)\,$. In
particular, this may happen if $\,\Omega\in\Lip_\alpha\,$, that
is, if $\,\partial\Omega\,$ coincides with the subgraph of a
$\,\Lip_\alpha$-function in a neighbourhood of each boundary
point. We prove that
$\,N_\NR(\Omega,\lambda)-C_{d,W}\,\mu_d(\Omega)\,\lambda^d=
O(\lambda^{(d-1)/\alpha})\,$ whenever $\,\Omega\in\Lip_\alpha\,$
and $\,\alpha\in(0,1)\,$ (Corollary \ref{C1.6}) and that this
estimate is order sharp (Theorem \ref{T1.10}). If
$\,(d-1)/\alpha<d\,$ then the right hand side is
$\,o(\lambda^d)\,$ and we have (\ref{0.1}), otherwise
$\,N_\NR(\Omega,\lambda)=O(\lambda^{(d-1)/\alpha})\,$.

We also obtain a remainder estimate in (\ref{0.1}) for the
Dirichlet Laplacian (Theorem \ref{T1.8}).  This estimate holds
true for all bounded domains and immediately implies (\ref{0.2}).

For domains with smooth boundaries our variational method only
gives the remainder estimate $\,O(\lambda^{d-1}\log\lambda)\,$; in
order to obtain $\,O(\lambda^{d-1})\,$ one has to use more
sophisticated results (see above). On the other hand, it can be
applied to many other problems and combined with the technique
developed in \cite{BI}, \cite{Iv3}, \cite{Iv4}, \cite{Me},
\cite{Mi}, \cite{SV} or \cite{Z} (see Section 5).

\bigskip
\noindent{\bf Acknowledgements.} The authors are very grateful to
M. Solomyak and E.B. Davies for their valuable comments.

\section{Definitions and main results}\label{S1}

\subsection{Basic definitions and notation}\label{S1.1}
Throughout the paper we assume that $\Omega$ is a bounded open
connected subset (domain) of the $d$-dimensional Euclidean space
$\R^d$ and that $d\geq2$.

We shall be using the following notation.
\begin{enumerate}
\item[$\bullet$]
$\omega_d$ is the volume of the unit ball in $\R^d$ and
$C_{d,W}:=(2\pi)^{-d}\,\omega_d$ is the standard Weyl constant.
\item[$\bullet$]
If $x=(x_1,\ldots,x_d)\in\R^d$ then $x':=(x_1,\ldots,x_{d-1})$ so
that $x=(x',x_d)$.
\item[$\bullet$]
$\,\overline\Omega\,$ and $\,\partial\Omega\,$ are the closure and
the boundary of $\Omega$.
\item[$\bullet$]
$\mu_d(\Omega)$ denotes the $d$-dimensional volume of $\Omega$ and
$D_\Omega:=\diam\Omega\,$.
\item[$\bullet$]
$\dist(\Omega_1,\Omega_2):=\inf\limits_{x\in\Omega_1,\,
y\in\Omega_2}|x-y|$ is the standard Euclidean distance.
\item[$\bullet$]
$\,\Omega_\eps^\bR:=\{x\in\Omega\;|\;\dist(x,\partial\Omega)\leq\eps\}\,$.
\item[$\bullet$]
$\,C\,$ is the space of continuous functions.
\item[$\bullet$]
If $\,\Omega'\,$ is a $\,(d-1)$-dimensional domain, $f\in
C(\overline{\Omega'})$, $b\in\R$ and $\,\alpha\in(0,1]\,$ then
\begin{align*}
\Gamma_f&:=\{x\in\R^d\;|\; x_d=f(x'),\ x'\in\Omega'\}\,,\\
G_f&:=\{x\in\R^d\;|\; x_d<f(x'),\ x'\in\Omega'\}\,,\\
G_{f,\,b}&:=\{x\in G_f\;|\; x_d>b\}\,,
\end{align*}
$\,\Osc(f,\Omega'):=\frac12\,\bigl(\sup\limits_{x\in\Omega'}f(x)
-\inf\limits_{x\in\Omega'} f(x)\bigr)\,$ and $\,|f|_\alpha\ :=\
\sup\limits_{x,\,y\in\Omega'}\frac{|f(x)-f(y)|}{|x-y|^\alpha}\,$.
\item[$\bullet$]
$\,Q_a^{(n)}\,$ is the open $\,n$-dimensional cube with edges of
length $a$ parallel to the coordinate axes. If the size or the
dimension of the cube $Q_a^{(n)}$ is not important for our
purposes or evident from the context then we shall omit the
corresponding index $a$ or $n$. However, we shall always be
assuming that the cube is open and that its edges are parallel to
the coordinate axes.
\item[$\bullet$]
$\Lip_\alpha$ is the space of functions $f$ on a cube $Q$ such
that $\,|f|_\alpha<\infty\,$ and $\lip_\alpha$ is the closure of
$\,\Lip_1\,$ in $\,\Lip_\alpha\,$ with respect to the seminorm
$\,|\cdot|_\alpha\,$.
\end{enumerate}

\begin{definition}\label{D1.1}
Given a bounded function $f$ on the cube $Q^{(n)}$ and $\delta>0$,
we shall denote by $\,\VC_\delta(f,Q^{(n)})\,$ the maximal number
of disjoint cubes $Q^{(n)}(i)\subset Q^{(n)}$ such that
$\Osc(f,Q^{(n)}(i))\geq \delta$ for each $i$. If
$\Osc(f,Q^{(n)})<\delta$ then we define
$\,\VC_\delta(f,Q^{(n)}):=1\,$.
\end{definition}

\begin{definition}\label{D1.2}
If $\,\tau\,$ is a positive nondecreasing function on
$\,(0,+\infty)\,$, let $\,BV_{\tau,\infty}(Q)\,$ be the space
spanned by all continuous functions $\,f\,$ on $\,\overline{Q}\,$
such that $\,\VC_{1/t}(f,Q)\leq\tau(t)\,$ for all $t>0\,$.
\end{definition}

We shall briefly discuss the relation between
$\,BV_{\tau,\infty}(Q)$ and known function spaces in
Subsection~\ref{S5.3}.

Let $\,X\,$ be a space of continuous real-valued functions defined
on a cube $\,Q^{(d-1)}\,$. We shall say that $\Omega$ belongs to
the class $X$ and write $\Omega\in X$ if for each
$z\in\partial\Omega$ there exists a neighbourhood $\,\OC_z\,$ of
the point $\,z\,$, a linear orthogonal map $U:\R^d\to\R^d$, a cube
$\,Q_a^{(d-1)}\subset Q^{(d-1)}\,$, a function $\,f\in X\,$ and
$\,b\in\R\,$ such that $\,U(\OC_z\bigcap\Omega)=\{x\in
G_{f,\,b}\;|\; x'\in Q_a^{(d-1)}\}\,$.

Since $\,\partial\Omega\,$ is compact, for every bounded set
$\Omega\in BV_{\tau,\infty}$ there exists a finite collection of
domains $\,\Omega_l\subset\Omega\,$, $\,l\in\LC\,$, such that
\begin{enumerate}
\item[(a)]
$\partial\Omega\subset\bigcup_{l\in\LC}\overline{\Omega_l}\,$;
\item[(b)]
for each $l$ we have $\,U_l(\Omega_l)=G_{f_l,\,b_l}\,$, where
$\,U_l:\R^d\to\R^d\,$ is a linear orthogonal map, $\,f_l\in
BV_{\tau,\infty}(Q_{a_l}^{(d-1)})\,$ and $\,b_l<\inf f_l\,$;
\item[(c)]
$\,a_l\leq D_\Omega\,$ and $\,\sup f_l-b_l\leq
D_\Omega\,$ for all $l\in\LC$.
\end{enumerate}
Let us fix such a collection $\,\{\Omega_l\}_{l\in\LC}\,$ and
denote $\,n_\Omega:=\#\LC\,$ and
$$
C_{\Omega,\,\tau}\ :=\ \sum_{l\in\LC}\sup_{t>0}
\left(\VC_{1/t}(f_l,Q_{a_l}^{(d-1)})/\tau(t)\right)\,.
$$
Let $\,\delta_\Omega\,$ be the largest positive number such that
$\,\Omega_{\delta_\Omega}^\bR\subset\bigcup_{l\in\LC}\Omega_l\,$,
$\,\delta_\Omega\leq\sqrt{d}\,a_l\,$ and $\,2\delta_\Omega\leq\inf
f_l-b_l\,$ for all $l\in\LC$.

\subsection{Main results}\label{S1.2} Throughout the paper we
shall denote by $C_d$ various constants depending only on the
dimension $d$. Constants appearing in the most important estimates
are numbered by an additional lower index; in our opinion, this
makes our proofs more transparent. Their precise (but not
necessarily best possible) values are given in Section 6.

\begin{theorem}\label{T1.3}
If $\,\Omega\in BV_{\tau,\infty}\,$ and
$\,\lambda\geq\delta_\Omega^{-1}\,$ then
\begin{multline}\label{1.1}
|\,N_\NR(\Omega,\lambda)-C_{d,W}\,\mu_d(\Omega)\,\lambda^d\,|\\
\leq\ C_{d,9}\,C_{\Omega,\tau}\,n_\Omega^{1/2}\lambda
\int_{(2D_\Omega)^{-1}}^{C_\Omega\,\lambda}t^{-2}\,\tau(t)\,\dR
t\;+\;C_{d,10}\,n_\Omega\,\lambda^{d-1} \int_0^{C_\Omega\,\lambda}
\mu_d(\Omega_{t^{-1}}^\bR)\,\dR t\,,
\end{multline}
where $\,C_\Omega:=4\,C_{d,8}\,n_\Omega^{1/2}\,$. If, in addition,
$\,\Omega\subset\R^2\,$ then there exists a positive constant
$\,c\,$ independent of $\Omega$ such that
\begin{multline}\label{1.2}
|\,N_\NR(\Omega,\lambda)-(4\pi)^{-1}\mu_2(\Omega)\,\lambda^2\,|\
\leq\ c\,C_{\Omega,\tau}\,\tau(c\,n_\Omega^{1/2}\lambda)\\
+\;c\,n_\Omega\,\lambda\,\left(D_\Omega+
\int_0^{c\,n_\Omega^{1/2}\lambda} \mu_2(\Omega_{t^{-1}}^\bR)\,\dR
t\right),\qquad\forall\lambda\geq \delta_\Omega^{-1}\,.
\end{multline}
\end{theorem}

\begin{remark}\label{R1.4}
For each continuous function $\,f\,$ on a closed cube there exists
a positive nondecreasing function $\,\tau\,$ such that $\,f\in
BV_{\tau,\infty}\,$. Therefore Theorem~\ref{T1.3} allows one to
obtain an estimate of the form (\ref{1.1}) for every domain
$\,\Omega\in C\,$. In particular, this implies the following well
known result: if $\,\Omega\in C\,$ then the essential spectrum of
the Neumann Laplacian on $\,\Omega\,$ is empty.
\end{remark}

The next two corollaries are simple consequences of
Theorem~\ref{T1.3}.

\begin{corollary}\label{C1.5}
If $\,\Omega\in BV_{\tau,\infty}\,$ then there exists a constant
$C_\Omega$ such that
\begin{multline}\label{1.3}
|\,N_\NR(\Omega,\lambda)-C_{d,W}\,\mu_d(\Omega)\,\lambda^d\,|\\
\leq\ C_\Omega\,\lambda^{d-1}
\int_{C_\Omega^{-1}}^{C_\Omega\lambda}\left(t^{-1}+t^{-d}\,\tau(t)\right)\dR
t\,, \qquad\forall\lambda\geq C_\Omega\,.
\end{multline}
\end{corollary}

\begin{corollary}\label{C1.6}
If $\,\alpha\in(0,1)$ and $\Omega\in\Lip_\alpha$ then
\begin{equation}\label{1.4}
N_\NR(\Omega,\lambda)\ =\ C_{d,W}\,\mu_d(\Omega)\,\lambda^d
\;+\;O\left(\lambda^{(d-1)/\alpha}\right),\qquad\lambda\to+\infty.
\end{equation}
If $\,\alpha\in(0,1)$ and $\Omega\in\lip_\alpha$ then
\begin{equation}\label{1.5}
N_\NR(\Omega,\lambda)\ =\ C_{d,W}\,\mu_d(\Omega)\,\lambda^d
\;+\;o\left(\lambda^{(d-1)/\alpha}\right),\qquad\lambda\to+\infty.
\end{equation}
\end{corollary}

\begin{remark}\label{R1.7}
If $\,\alpha\leq1-d^{-1}\,$ then the asymptotic formula
(\ref{1.4}) turns into the estimate
$\,N_\NR(\Omega,\lambda)=O\left(\lambda^{(d-1)/\alpha}\right)$.
Similarly, if $\,\alpha<1-d^{-1}\,$ then (\ref{1.5}) takes the
form
$\,N_\NR(\Omega,\lambda)=o\left(\lambda^{(d-1)/\alpha}\right)$.
\end{remark}

The following estimates for the Dirichlet Laplacian are much
simpler. The inequality (\ref{1.6}) seems to be new but results of
this type are known to experts. Corollary~\ref{C1.9} is an
immediate consequence of Theorem \ref{T1.8}; (\ref{1.7}) also
follows from (\ref{0.2}).

\begin{theorem}\label{T1.8}
For all $\,\lambda>0\,$ we have
\begin{equation}\label{1.6}
|\,N_\DR(\Omega,\lambda)-C_{d,W}\,\mu_d(\Omega)\,\lambda^d\,|\
\leq\ C_{d,11}\,\lambda^{d-1}
\int_0^\lambda\mu_d(\Omega_{t^{-1}}^\bR)\,\dR t\,.
\end{equation}
\end{theorem}

\begin{corollary}\label{C1.9}
If $\,\alpha\in(0,1)$ and $\Omega\in\Lip_\alpha$ then
\begin{equation}\label{1.7}
N_\DR(\Omega,\lambda)\ =\ C_{d,W}\,\mu_d(\Omega)\,\lambda^d
\;+\;O\left(\lambda^{d-\alpha}\right),\qquad\lambda\to+\infty.
\end{equation}
If $\,\alpha\in(0,1)$ and $\Omega\in\lip_\alpha$ then
\begin{equation}\label{1.8}
N_\DR(\Omega,\lambda)\ =\ C_{d,W}\,\mu_d(\Omega)\,\lambda^d
\;+\;o\left(\lambda^{d-\alpha}\right),\qquad\lambda\to+\infty.
\end{equation}
\end{corollary}

Note that $\,(d-1)/\alpha>d-\alpha\,$ whenever
$\,\alpha\in(0,1)\,$. Therefore the remainder estimate in
Corollary~\ref{C1.9} is better than that in Corollary~\ref{C1.6}.
The following theorem shows that the asymptotic formulae
(\ref{1.4}) and (\ref{1.5}) are order sharp.

\begin{theorem}\label{T1.10}
Let $\alpha\in(0,1)$. Then
\begin{enumerate}
\item[(1)] there exist a bounded domain $\Omega\in\Lip_\alpha$ and
a positive constant $C_\Omega$ such that
$\;N_\NR(\Omega,\lambda)\geq
C_{d,W}\,\mu_d(\Omega)\,\lambda^d+C_\Omega^{-1}\,\lambda^{(d-1)/\alpha}\;$
for all $\,\lambda>C_\Omega\,$;
\item[(2)]
for each nonnegative function $\phi$ on $(0,+\infty)$ vanishing at
$+\infty$ there exist a bounded domain $\Omega\in\lip_\alpha$ and
a positive constant $C_{\phi,\Omega}$ such that
$\;N_\NR(\Omega,\lambda)\geq C_{d,W}\,\mu_d(\Omega)\,\lambda^d+
C_{\phi,\Omega}^{-1}\,\phi(\lambda)\,\lambda^{(d-1)/\alpha}\;$ for
all $\,\lambda>C_{\phi,\Omega}\,$.
\end{enumerate}
\end{theorem}

\begin{remark}\label{R1.11}
In \cite{BD} the authors proved that
\begin{equation}\label{1.9}
0<K_{\Omega,\NR}(t,x,y)\leq
C_\Omega\,t^{-(\alpha+d-1)/(2\alpha)}\,,\qquad\forall
x,y\in\Omega,\quad\forall t\in(0,1],
\end{equation}
whenever $\,\Omega\in\Lip_\alpha\,$ and $\,\alpha\in(0,1)\,$,
where $\,K_{\Omega,\NR}\,$ is the heat kernel of the Neumann
Laplacian on $\,\Omega\,$ and $\,C_\Omega\,$ is a constant
depending on $\,\Omega\,$. The estimate (\ref{1.9}) is order sharp
as $\,t\to0\,$ (see \cite{BD}, Example 6). Corollary~\ref{C1.6}
implies that there exists a constant $\,C'_\Omega\,$ such that
$$
\int_\Omega K_{\Omega,\NR}(t,x,x)\,\dR x\leq
C'_\Omega\,(t^{-d/2}+t^{-(d-1)/(2\alpha)})\,,\qquad\forall
t\in(0,1].
$$
In view of Theorem~\ref{T1.10}, this estimate is also order sharp.
Since $\,d/2<(\alpha+d-1)/(2\alpha)\,$ and
$\,(d-1)/(2\alpha)<(\alpha+d-1)/(2\alpha)\,$, we see that
integration of the heat kernel $\,K_{\Omega,\NR}(t,x,x)\,$
improves its asymptotic properties.
\end{remark}

\subsection{Further definitions and notation}\label{S1.3}
In the rest of the paper
\begin{enumerate}
\item[$\bullet$]
$\#T$ denotes the number of elements of the set $T$.
\item[$\bullet$]
If $\,\{T(i)\}_{i\in\IC}\,$ is a finite family of sets $\,T(i)\,$
and $\,T:=\bigcup_{i\in\IC} T(i)\,$ then
$$
\aleph\{T(i)\}:=\sup_{x\in T}\left(\#\{i\in\IC\,|\,x\in
T(i)\}\right),
$$
in other words, $\,\aleph\{T(i)\}\,$ is the multiplicity of the
covering $\,\{T(i)\}_{i\in\IC}\,$.
\item[$\bullet$]
If $\,s\in\R_+\,$ then $\,[s]\,$ is the entire part of $\,s\,$.
\item[$\bullet$]
$\supp f$ and $\nabla f$ denote the support and gradient of the
function $f$.
\end{enumerate}

The paper is organised as follows. In the next section we recall
some well known results from spectral theory and estimate the
counting function on `model' domains. In Section 3 we discuss
partitions of the domain $\Omega\,$. In Section 4 we deduce the
main theorems from the results of Sections 2 and 3. In the last
section we extend our results to a wider class of domains and
higher order operators and discuss other possible generalizations.

\section{Variational formulae and related results}\label{S2}

Recall that the Sobolev space $W^{1,2}(\Omega)$ is the space of
functions $u\in L^2(\Omega)$ such that $\nabla u\in L^2(\Omega)$,
endowed with the norm
$$
\|u\|_{W^{1,2}(\Omega)}=(\|\nabla u\|_{L_2(\Omega)}^2
+\|u\|_{L^2(\Omega)}^2)^{1/2}.
$$
If $\Upsilon$ is a subset of $\partial\Omega$, let
$W^{1,2}_{0,\Upsilon}(\Omega)$ be the closure in $W^{1,2}(\Omega)$
of the set
$$
\{f\in W^{1,2}(\Omega)\;|\;\;\supp f\bigcap \Upsilon=\emptyset\}
$$
and $W_0^{1,2}(\Omega):=W^{1,2}_{0,\partial\Omega}(\Omega)$.
Obviously, $W^{1,2}_{0,\emptyset}(\Omega)=W^{1,2}(\Omega)$.

Let
\begin{equation}\label{2.1}
N_{\NR,\DR}(\Omega,\Upsilon,\lambda)\ :=\ \sup(\dim E_\lambda)
\end{equation}
where the supremum is taken over all subspaces $E_\lambda\subset
W^{1,2}_{0,\Upsilon}(\Omega)$ such that
\begin{equation}\label{2.2}
\|\nabla u\|^2_{L_2(\Omega)}\ <\
\lambda^2\,\|u\|^2_{L^2(\Omega)}\,,\qquad\forall u\in E_\lambda\,.
\end{equation}
In view of the Rayleigh--Ritz variational formula,
$N_{\NR,\DR}(\Omega,\Upsilon,\lambda)$ can be thought of as the
counting function of the Laplacian on the bounded domain $\Omega$
subject to Dirichlet boundary condition on $\Upsilon$ and Neumann
boundary condition on the remaining part of the boundary. In
particular,
$\,N_{\NR,\DR}(\Omega,\emptyset,\lambda)=N_\NR(\Omega,\lambda)\,$
and
$\,N_{\NR,\DR}(\Omega,\partial\Omega,\lambda)=N_\DR(\Omega,\lambda)\,$.
Equivalently, (\ref{2.1}) can be rewritten as
\begin{equation}\label{2.3}
N_{\NR,\DR}(\Omega,\Upsilon,\lambda)\ =\ \inf(\codim\tilde
E_\lambda),
\end{equation}
where the infimum is taken over all subspaces $\tilde
E_\lambda\subset W^{1,2}_{0,\Upsilon}(\Omega)$ such that
\begin{equation}\label{2.4}
\|\nabla u\|^2_{L_2(\Omega)}\ \geq\
\lambda^2\,\|u\|^2_{L^2(\Omega)}\,,\qquad\forall u\in\tilde
E_\lambda\,.
\end{equation}

\begin{lemma}\label{L2.1}
Let $\{\Omega_i\}_{i\in\IC}$ be a countable family of disjoint
open sets $\Omega_j\subset\Omega$ such that
$\mu_d(\Omega)=\mu_d(\bigcup_{i\in\IC}\Omega_i)$. Then
$$
\sum_{i\in\IC}N_\DR(\Omega_i,\lambda)\ \leq\
N_\DR(\Omega,\lambda)\ \leq\ N_\NR(\Omega,\lambda)\ \leq\
\sum_{i\in\IC}N_\NR(\Omega_i,\lambda)
$$
and $\,N_\NR(\Omega,\lambda)\geq\sum_{j\in\JC}N_{\NR,\DR}
(\Omega_j,\partial\Omega_j\setminus\partial\Omega,\lambda)\,$.
\end{lemma}

Lemma~\ref{L2.1} is an elementary corollary of the Rayleigh--Ritz
formula. The following lemma is less obvious.

\begin{lemma}\label{L2.2}
Let $\{\Omega_i\}_{i\in\IC}$ be a countable family of open sets
$\Omega_j\subset\Omega$ such that
$\,\mu_d(\Omega)=\mu_d(\bigcup_{i\in\IC}\Omega_i)\,$,
$\,\Upsilon\,$ be an arbitrary subset of $\,\partial\Omega\,$ and
$\,\Upsilon_j:=\partial\Omega_j\bigcap\Upsilon\,$. If
$\,\aleph\{\Omega_j\}\leq\varkappa<+\infty$ then
$N_{\NR,\DR}(\Omega,\Upsilon,\varkappa^{-1/2}\lambda)
\leq\sum_{j\in\JC}N_{\NR,\DR}(\Omega_j,\Upsilon_j,\lambda)$.
\end{lemma}

\begin{proof}
Denote by $\,\tilde E_{\lambda,j,\Omega}\,$ the subspace of
functions $\,u\in W_{0,\Upsilon}^{1,2}(\Omega)\,$ such that
$\,\|\nabla u\|^2_{L_2(\Omega_j)}\geq
\lambda^2\,\|u\|^2_{L^2(\Omega_j)}\,$. We have
$\,\varkappa\,\|u\|^2_{L^{1,2}(\Omega)}
\geq\lambda^2\,\|u\|^2_{L^2(\Omega)}\,$ whenever
$u\in\bigcap_{j\in\JC}\tilde E_{\lambda,j,\Omega}\,$. Therefore,
by (\ref{2.3}),
$$
N_\NR(\Omega,\varkappa^{-1/2}\lambda)\ \le\
\inf(\codim\bigcap_{j\in\JC}\tilde E_{\lambda,j,\Omega})\ \le\
\sum_{j\in\JC}\inf(\codim\tilde E_{\lambda,j,\Omega}),
$$
where the infimum are taken over all subspaces $\tilde
E_{\lambda,j,\Omega}$ satisfying the above condition. If $\tilde
E_{\lambda,j}$ is the intersection of the kernels of linear
continuous functionals $\Lambda_k$ on
$W_{0,\Upsilon_j}^{1,2}(\Omega_k)$ and $E_{\lambda,j,\Omega}$ is
the intersection of the kernels of linear continuous functionals
$u\to\Lambda_k(\left.u\right|_{\Omega_j})$ on
$W_{0,\Upsilon}^{1,2}(\Omega)$ then $\codim\tilde
E_{\lambda,j}\geq\codim E_{\lambda,j,\Omega}$ and
$\left.u\right|_{\Omega_j}\in\tilde E_{\lambda,j}$ whenever $u\in
E_{\lambda,j,\Omega}$. This observation and (\ref{2.3}) imply that
$\inf(\codim\tilde E_{\lambda,j,\Omega})\le
N_{\NR,\DR}(\Omega_j,\Upsilon_j,\lambda)$.
\end{proof}

\begin{remark}\label{R2.3}
Lemma~\ref{L2.2} implies that
$\,N_\NR(\Omega,\varkappa^{-1/2}\lambda)
\leq\sum_{j\in\JC}N_\NR(\Omega_j,\lambda)\,$ whenever
$\,\bigcup_{j\in\JC}\Omega_j\subset\Omega$,
$\mu_d(\Omega)=\mu_d(\bigcup_{i\in\IC}\Omega_i)$ and
$\,\aleph\{\Omega_j\}\leq\varkappa\,$. It may well be the case
that, under these conditions,
$N_\NR(\Omega,\lambda)\leq\sum_{j\in\JC}N_\NR(\Omega_j,\lambda)\,$.
This conjecture looks plausible and is equivalent to the following
statement: if $\Omega_1\subset\Omega$, $\Omega_2\subset\Omega$ and
$\mu_d(\Omega)=\mu_d(\Omega_1)+\mu_d(\Omega_2)$ then
$N_\NR(\Omega_1,\lambda)+N_\NR(\Omega_2,\lambda)\geq
N_\NR(\Omega,\lambda)$.
\end{remark}

\begin{remark}\label{R2.4}
The first eigenvalue of the Neumann Laplacian $-\Delta _\NR$ is
always equal to 0 and the corresponding eigenfunction is
identically equal to constant. Let $\,\lambda_{1,\NR}(\Omega)
:=\inf\{\lambda\in\R_+\,|\,N_\NR(\Omega,\lambda)>1\}\,$; if
$-\Delta _\NR$ has at least two eigenvalues lying below its
essential spectrum (or the essential spectrum is empty) then
$\,\lambda_{1,\NR}(\Omega)\,$ coincides with the smallest nonzero
eigenvalue of the operator $\,\sqrt{-\Delta _\NR}\,$. By the
spectral theorem, we have $\,\lambda_{1,\NR}(\Omega)\geq\lambda\,$
if and only if $\,\int_\Omega|u(x)|^2\,\dR x\leq
\lambda^{-2}\int_\Omega|\nabla u(x)|^2\,\dR x\,$ for all functions
$\,u\in W^{1,2}(\Omega)\,$ such that $\,\int_\Omega u(x)\,\dR
x=0\,$. Note that $\,\int_\Omega|u(x)|^2\,\dR
x\leq\int_\Omega|u(x)-c|^2\,\dR x\,$ for all $c\in\C$ whenever
$\int_\Omega u(x)\,\dR x=0$.
\end{remark}

\begin{definition}\label{D2.5}
Denote by $\PB(\delta)$ the set of all rectangles with edges
parallel to the coordinate axes, such that the length of the
maximal edge does not exceed $\delta\,$. If $f$ is a continuous
function on $\overline{Q^{(d-1)}}$, let $\VB(\delta,f)$ be the
class of domains $V\subset G_f$ which can be represented in the
form $\,V=G_{f,\,b}(Q_c^{(d-1)})\,$, where $\,Q_c^{(d-1)}\subset
Q^{(d-1)}\,$, $\,c\leq\delta\,$, $\,b=\inf f-\delta\,$ and
$\,\Osc(f,Q_c^{(d-1)})\leq\delta/2\,$. We shall write
$\,V\in\VB(\delta)\,$ if $\,V\in\VB(\delta,f)\,$ for some
continuous function $f$. Finally, let $\,\MB(\delta)\,$ be the
class of open sets $\,M\subset\R^d\,$ such that $\,M\subset
Q_\delta^{(d)}\,$ for some cube $\,Q_\delta^{(d)}\,$.
\end{definition}

\begin{lemma}\label{L2.6} Let $\,\delta\,$ be an arbitrary
positive number.
\begin{enumerate}
\item[(1)]
If $\,P\in\PB(\delta)$ then $N_\NR(P,\lambda)=1\,$ for all
$\,\lambda\leq\pi\delta^{-1}$.
\item[(2)]
If $V\in\VB(\delta)$ then $N_\NR(V,\lambda)=1\,$ for all
$\lambda\leq(1+2\pi^{-2})^{-1/2}\delta^{-1}$.
\item[(3)]
If $\,M\in\MB(\delta)\,$, $\,M\subset Q_\delta^{(d)}\,$ and
$\,\Upsilon:=\partial M\bigcap Q_\delta^{(d)}\,$ then we have
$\,N_{\NR,\DR}(M,\Upsilon,\lambda)\leq1\,$ for all
$\,\lambda\leq\pi\delta^{-1}$ and
$\,N_{\NR,\DR}(M,\Upsilon,\lambda)=0\,$ for all $\,\lambda\leq
(2^{-1}-2^{-1}\delta^{-d}\mu_d(M))^{1/2}\,\pi\delta^{-1}\,$.
\end{enumerate}
\end{lemma}

\begin{proof}
If $\,P\,$ is a rectangle then $\,\lambda_{1,\NR}=\pi\,a^{-1}\,$,
where $\,a\,$ is the length of its maximal edge. This implies (1).

Assume now that $\,V\in\VB(\delta,f)\,$, where $\,f\,$ is a
continuous function on $\,\overline{Q_c^{(d-1)}}\,$ and denote
$\,b:=\inf f-\delta\,$ and $\,P:=Q_c^{(d-1)}\times(b,b+\delta)\,$.
Clearly, $\,P\in\PB(\delta)\,$. Let $\,u\in W^{1,2}(V)\,$ and
$\,c'_u\,$ the average of $\,u\,$ over $\,P\,$. If
$\,r\in[b,b+\delta]\,$ and $\,s\in[b+\delta,f(x')]\,$ then, by
Jensen's inequality,
$$
|u(x',s)-u(x',r)|^2\ =\ |\int_r^s\partial_t\,u(x',t)\,\dR t\,|^2\
\leq\ (s-r)\int_b^{f(x')}|\partial_t\,u(x',t)|^2\,\dR t\,.
$$
Since $\,\int_b^{b+\delta}\int_{b+\delta}^f(s-r)\,\dR s\,\dR r=
(\delta/2)\,(f-b-\delta)\,(f-b)\,$ and
$$
0\leq f-b-\delta=f-\inf f\leq2\,\Osc(f,Q_c^{(d-1)})\leq\delta\,,
$$
we have
$$
\int_b^{g(x')}\int_{g(x')}^{f(x')}|u(x',s)-u(x',r)|^2\,\dR s\,\dR r\\
\leq\ \delta^3\int_b^{f(x')}|\partial_t\,u(x',t)|^2\,\dR t\,.
$$
In view of Remark~\ref{R2.4} and (1), we also have
\begin{equation}\label{2.5}
\int_P|u(x)-c'_u|^2\,\dR x\ \leq\ \pi^{-2}\,\delta^2\int_P|\nabla
u(x)|^2\,\dR x.
\end{equation}
Integrating the inequality
$$
|u(x',s)-c'_u|^2\ \leq\
(1+\gamma)\,|u(x',r)-c'_u|^2+(1+\gamma^{-1})\,|u(x',s)-u(x',r)|^2
$$
over $\,r\in[b,b+\delta]\,$, $\,s\in[b+\delta,f(x')]\,$ and
$\,x'\in\Omega'\,$ and applying these two estimates, we obtain
\begin{multline*}
\delta\int_{V\setminus P}|u(x)-c'_u|^2\,\dR x\ \leq
(1+\gamma)\,\pi^{-2}\,\delta^3
\int_P|\nabla u(x)|^2\,\dR x\\
+(1+\gamma^{-1})\,\delta^3 \int_V|\partial_{x_d}u(x)|^2\,\dR x
\end{multline*}
for all $\,\gamma>0\,$. Dividing both sides by $\,\delta\,$ and
substituting $\,\gamma=\pi^2\,$, we see that $\,\int_{V\setminus
P}|u(x)-c'_u|^2\,\dR x\,$ is estimated by
$\,(1+\pi^{-2})\,\delta^2\int_V|\nabla u(x)|^2\,\dR x\,$. Now (2)
follows from (\ref{2.5}) and Remark~\ref{R2.4}.

In order to prove (3), let us consider a function $\,u\in
W^{1,2}(M)\,$ which vanishes near $\,\Upsilon\,$ and extend it by
zero to the whole cube $\,Q_\delta^{(d)}\,$. Since $\,u\in
W^{1,2}(Q_\delta^{(d)})\,$, (1) implies the first inequality (3).
If $\,c_u\,$ is the average of $\,u\,$ over $\,Q_\delta^{(d)}\,$
then
\begin{equation}\label{2.6}
\int_M|c_u|^2\,\dR x\ \leq\
\mu_d(M)\,\delta^{-d}\left(\int_M|c_u|^2\,\dR x
+\int_{Q_\delta^{(d)}}|u(x)-c_u|^2\,\dR x\right).
\end{equation}
Therefore Remark~\ref{R2.4} and (1) imply that
\begin{multline*}
\int_M|u(x)|^2\,\dR x\ \leq\
2\int_{Q_\delta^{(d)}}|u(x)-c_u|^2\,\dR
x+2\int_M|c_u|^2\,\dR x\\
\leq\ 2\left(1+\mu_d(M)\,\delta^{-d}
\left(1-\mu_d(M)\,\delta^{-d}\right)^{-1}\right)\int_{Q_\delta^{(d)}}
|u(x)-c_u|^2\,\dR x\\ \leq\
2\,\pi^{-2}\,\delta^2\left(1-\mu_d(M)\,\delta^{-d}\right)^{-1}
\int_M|\nabla u(x)|^2\,\dR x\,.
\end{multline*}
The second identity (3) follows from the above inequality and the
Rayleigh--Ritz formula.
\end{proof}

\begin{remark}\label{R2.7}
The second estimate in Lemma~\ref{L2.6}(3) is sufficient for our
purposes but is very rough. One can obtain a much more precise
result in terms of capacities (see \cite{M2}, Chapter 10, Section
1).
\end{remark}

\begin{lemma}\label{L2.8}
Let $\,\delta>0\,$. Then for all $\,\lambda>0\,$ we have
$$
-\,C_{d,1}\left((\delta\lambda)^{d-1}+1\right)\ \leq\
N(Q_\delta^{(d)},\lambda)-C_{d,W}\,(\delta\lambda)^d\ \leq\
C_{d,1}\left((\delta\lambda)^{d-1}+1\right).
$$
\end{lemma}

\begin{proof}
Changing variables $\tilde x=\delta\,x$, we see that
\begin{equation}\label{2.7}
N(\Omega,\delta\lambda)\ =\
N(\delta\Omega,\lambda)\,,\quad\text{where}\quad
\delta\Omega:=\{x\in\R^d\,|\,\delta^{-1}x\in\Omega\}\,.
\end{equation}
Therefore it is sufficient to prove the required estimates only
for $\,\delta=1\,$. If $\,\Omega=\Omega'\times\Omega''\,$,
$\,\Upsilon'\subset\partial\Omega'\,$ and
$\,\Upsilon''\subset\partial\Omega''\,$ then, separating
variables, we obtain
\begin{equation}\label{2.8}
N_{\NR,\DR}(\Omega,\Upsilon,\lambda)\ =\ \int
N_{\NR,\DR}\left(\Omega',\Upsilon',\sqrt{\lambda^2-\mu^2}\right)\,\dR
N_{\NR,\DR}(\Omega'',\Upsilon'',\mu)\,,
\end{equation}
where $\,\Upsilon=(\Upsilon'\times\partial\Omega'')\bigcup
(\partial\Omega'\times\Upsilon'')\,$ and the right hand side is a
Stieltjes integral. Using (\ref{2.8}), explicit formulae for the
counting functions on the unit interval and the identities
\begin{equation}\label{2.9}
\int_0^\lambda(\lambda^2-\mu^2)^{n/2}\,\dR\mu \ =\
\lambda^{n+1}\,\omega_{n+1}\,(2\,\omega_n)^{-1}\,,\qquad \forall
n=1,2,\ldots,
\end{equation}
one can easily prove the required inequality by induction in
$\,d\,$.
\end{proof}

\begin{remark}\label{R2.9}

Lemma~\ref{L2.8} is an immediate consequence of well known results
on spectral asymptotics in domains with piecewise smooth
boundaries (see, for example, \cite{Iv2} or \cite{F}); a similar
result holds true for higher order elliptic operators and
operators with variable coefficients \cite{V}. We have given an
independent proof in order to find the explicit constant
$\,C_{d,1}\,$.
\end{remark}

\section{Properties of domains and their partitions}\label{S3}

\subsection{Besicovitch's and Whitney's theorems}\label{S3.1}
We shall use the following version of Besicovitch's theorem.

\begin{theorem}\label{T3.1}
There are two constants $\,\CC_n\geq1\,$ and $\,\hat\CC_n\geq1\,$
depending only on the dimension $n$, such that for every compact
set $\,K\subset\R^n\,$ and every positive function $\rho$ on $K$
one can find a finite subset $\,\YC\subset K\,$ and a family of
cubes $\,\{Q_{\rho(y)}^{(n)}[y]\}_{y\in\YC}\,$ centred on $y$,
which satisfy the following conditions:
\begin{enumerate}
\item[(1)]
$K\subset\bigcup_{y\in\YC}Q_{\rho(y)}^{(n)}[y]\,$,
\item[(2)]
$\aleph\{K\bigcap Q_{\rho(y)}^{(n)}[y]\}_{y\in\YC}\leq\CC_n\,$;
\item[(3)]
there exists a subset $\hat\YC\subset\YC$ such that $\,\#\YC\leq
\hat\CC_n(\#\hat\YC)\,$ and the cubes
$\{Q_{\rho(y)}^{(n)}[y]\}_{y\in\hat{\YC}}$ are mutually disjoint.
\end{enumerate}
\end{theorem}

Theorem \ref{T3.1} is proved in the same way as Besicovitch's
theorem in \cite{G}, Chapter 1.

\begin{corollary}\label{C3.2}
Let $f$ be a continuous function on the closure
$\overline{Q^{(d-1)}}$. Then for every $\eps>0$ there exists a
finite family of cubes $\{Q^{(d-1)}(x)\}_{x\in\XC}$ such that
\begin{enumerate}
\item[(1)]
$\bigcup_{x\in\XC}\overline{Q^{(d-1)}}(x)=\overline{Q^{(d-1)}}$;
\item[(2)]
$\aleph\{Q^{(d-1)}(x)\}\leq C_{d,2}$;
\item[(3)]
$\#\XC\leq C_{d,3}\,\VC_\eps(f,\overline{Q^{(d-1)}})$;
\item[(4)]
$\Osc(f,Q^{(d-1)}(x))\leq\eps$ for each $x\in\XC$.
\end{enumerate}
\end{corollary}

\begin{proof}
Without loss of generality we can assume that
$Q^{(d-1)}=(-1,1)^{d-1}$ and $\Osc(f,Q^{(d-1)})>\eps$. Let us
denote by $\,Q_t^{(d-1)}[y]\,$ the cube of the size $\,t\,$
centred on $\,y\,$, define
$$
\rho(y):=\inf\{t>0\;|\;\Osc(f,Q^{(d-1)}\bigcap
Q_t^{(d-1)}[y])=\eps\}\,,\qquad y\in\overline{Q^{(d-1)}}\,,
$$
apply Besicovitch's theorem to the set $K=\overline{Q^{(d-1)}}$
and find the sets $\YC$ and $\hat\YC$. If $y\in\YC$, denote
$\,P^{(d-1)}[y]:=Q^{(d-1)}\bigcap Q_{\rho(y)}^{(d-1)}[y]\,$ and
assume that
$$
P^{(d-1)}[y]\ =\ (a_1(y),b_1(y))\times(a_2(y),b_2(y))
\times\dots\times(a_{d-1}(y),b_{d-1}(y))\,,
$$
where $\,-1\leq a_j(y)<b_j(y)\leq1\,$. Let $Q'(y)$ be the minimal
cube such that $P^{(d-1)}(x)\subset Q'(y)\subset Q^{(d-1)}$ and
$c(y):=\max_j\,(b_j(y)-a_j(y))$. We have
$$
Q'(y)\ =\ (a'_1(y),b'_1(y))\times(a'_2(y),b'_2(y))
\times\dots\times(a'_{d-1}(y),b'_{d-1}(y))\,,
$$
where
\begin{enumerate}
\item[(-1)]
if $a_j(y)=-1$ then $a'_j(y)=-1$ and
$b'_j(y)=a_j(y)+c(y)$;
\item[(0)]
if $a_j(y)>-1$ and $b_j(y)<1$ then $a'_j(y)=a_j(y)$ and
$b'_j(y)=b_j(y)$;
\item[(+1)]
if $b_j(y)=1$ then $a'_j(y)=b_j(y)-c(y)$ and $b'_j(y)=1$.
\end{enumerate}

Let us consider the set $\Sigma=\{-1,0,1\}^{d-1}$ of all
$\,(d-1)$-dimensional vectors
$\sigma=(\sigma_1,\ldots,\sigma_{d-1})$ with entries $\sigma_j$
equal to $\,-1$, $0$ or $1$. Denote by $\hat\YC_\sigma$ the set of
points $y\in\hat\YC$ such that $a_j(y)$ and $b_j(y)$ satisfy the
condition ($\sigma_j$) for all $j=1,\ldots,d-1$. Since
$\,\aleph\{P^{(d-1)}[y]\}_{y\in\hat\YC}=1$, for each
$\sigma\in\Sigma$ the cubes $\{Q'(y)\}_{y\in\hat\YC_\sigma}=1$ are
mutually disjoint. Therefore $\#\hat\YC_\sigma\leq
\VC_\eps(f,\overline{Q^{(d-1)}})$ for all $\sigma\in\Sigma$ (see
Definition~\ref{D1.1}) and, consequently,
$\,\#\hat\YC\leq(\#\Sigma)\,\VC_\eps(f,\overline{Q^{(d-1)}})
\leq3^{d-1}\,\VC_\eps(f,\overline{Q^{(d-1)}})\,$. This estimate
and Theorem~\ref{T3.1}(3) imply that
$\,\#\YC\leq3^{d-1}\,\hat\CC_{d-1}\,\VC_\eps(f,\overline{Q^{(d-1)}})\,$.

Since $\,\YC\subset\overline{Q^{(d-1)}}\,$, we have
$\,1/2\leq(b_j(y)-a_j(y))^{-1}(b_k(y)-a_k(y))\leq2\,$ for all
$\,j,k=1,\ldots,d-1\,$ and $\,y\in\YC\,$. Using this inequality,
one can easily show by induction in $\,d\,$ that every rectangle
$\,P^{(d-1)}[y]\,$ coincides with the union of a finite collection
of cubes $\,\{Q^{(d-1)}(x)\}_{x\in\XC_y}\,$ such that
$\,\#\XC_y\leq2^{d-1}\,$ and
$\,\aleph\{Q^{(d-1)}(x)\}_{x\in\XC_y}\leq2^{d-1}\,$.

Let $\,\XC:=\bigcup_{y\in\YC}\XC_y\,$. In view of the first two
conditions of Theorem~\ref{T3.1}, the family
$\,\{Q^{(d-1)}(x)\}_{x\in\XC}\,$ satisfies (1) and (2). The upper
bound
$\,\#\YC\leq3^{d-1}\,\hat\CC_{d-1}\,\VC_\eps(f,\overline{Q^{(d-1)}})\,$
implies (3). Finally, since $\,\Osc(f,P^{(d-1)}[y])=\eps\,$ and
$\,Q^{(d-1)}(x)\subset P^{(d-1)}[y]\,$ whenever $\,x\in\XC_y\,$ ,
we have (4).
\end{proof}

The following theorem is due to Whitney. It can be found, for
example, in \cite{St}, Chapter VI, or \cite{G}, Chapter 1.

\begin{theorem}\label{T3.3}
There exists a countable family of mutually disjoint cubes
$\,\{Q_{2^{-i}}^{(d)}(i,n)\}_{n\in\NC(i)\,,\,i\in\IC}\,$ such that
$\,\overline{\Omega}=\bigcup_{i\in\IC}\bigcup_{n\in\NC_i}
\overline{Q_{2^{-i}}^{(d)}(i,n)}\,$ and
\begin{equation}\label{3.1}
Q_{2^{-i}}^{(d)}(i,n)\subset\{x\in\Omega\;|\;
\sqrt{d}\,2^{-i}\leq\dist(x,\partial\Omega)\leq4\sqrt{d}\,2^{-i}\}\,.
\end{equation}
Here $\IC$ is a subset of $\Z$ and $\NC_i$ are some finite index
sets.
\end{theorem}

\subsection{Auxiliary results}\label{S3.2}
In this subsection we shall prove several technical results
concerning domains $G_{f,\,b}\,$.

\begin{lemma}\label{L3.4}
Let $f$ be a continuous function defined on the closure
$\overline{Q_a^{(d-1)}}$. Then for every $\delta>0$ and $m\in\Z_+$
there exists a finite family of cubes
$\{Q^{(d-1)}(k)\}_{k\in\KC_m}$ such that
\begin{enumerate}
\item[(1)]
$\bigcup_{k\in\KC_m}\overline{Q^{(d-1)}(k)}=\overline{Q_a^{(d-1)}}$;
\item[(2)] $Q^{(d-1)}(k)\in\PB(\delta)$ for all $k\in\KC_m$;
\item[(3)] $\,\aleph\{Q^{(d-1)}(k)\}_{k\in\KC_m} \leq C_{d,2}$;
\item[(4)] $\Osc(f,Q^{(d-1)}(k))\leq2^{m-1}\delta$ for all $k\in
\KC_m$; \item[(5)]
$\#\{k\in\KC_m\;|\;\mu_{d-1}(Q^{(d-1)}(k))\leq2^{1-d}\,\delta^{d-1}\}\leq
C_{d,3}\,\VC_{2^{m-1}\delta}(f,Q_a^{(d-1)})\,$.
\end{enumerate}
\end{lemma}

\begin{proof}
Let $\{Q^{(d-1)}(x)\}_{x\in\XC}$ be a family of cubes satisfying
the conditions of Corollary~\ref{C3.2} with $\eps=2^{m-1}\delta$.
Assume that $Q^{(d-1)}(x)=Q_{a_x}^{(d-1)}\,$ with some $a_x>0$ and
denote by $\XC_\delta$ the set of all indices $x\in\XC$ such that
$a_x\leq\delta$. For each $x\in\XC\setminus\XC_\delta$, we choose
a positive integer $m_x$ such that $a_x/m_x\in(\delta/2,\delta]$
and split the closed cube $\overline{Q^{(d-1)}(x)}$  into the
union of $m_x^{d-1}$ congruent closed cubes
$\overline{Q_{a_x/m_x}^{(d-1)}(x,j)}$, $j=1,\ldots,m_x^{d-1}$. Let
$Q_{a_x/m_x}^{(d-1)}(x,j)$ be the corresponding disjoint open
cubes and
$$
\{Q^{(d-1)}(k)\}_{k\in\KC}:=\{Q^{(d-1)}(k)\}_{x\in\XC_\delta}
\bigcup\{Q_{a_x/m_x}^{(d-1)}(x,j)\}_{x\in\XC\setminus\XC_\delta,
\,j=1,\ldots,m_x^{d-1}}\,.
$$
Then (2) holds true and  (1), (3), (4) and (5) follow from
Corollary~\ref{C3.2}(1), Corollary~\ref{C3.2}(2),
Corollary~\ref{C3.2}(4) and Corollary~\ref{C3.2}(3) respectively.
\end{proof}

\begin{theorem}\label{T3.5}
Let $f$ be a continuous function on $\overline{Q_a^{(d-1)}}\,$,
$\,\delta\in(0,\sqrt{d}\,a]\,$ and $\,b\in[-\infty,\,\inf
f-2\delta]\,$. Then there exist countable families of sets
$\{P_j\}_{j\in\JC}$ and $\{V_k\}_{k\in\KC}$ satisfying the
following conditions:
\begin{enumerate}
\item[(1)] $P_j\subset G_{f,b}$ and $P_j\in\PB(\delta)$ for all
$j\in\JC$; \item[(2)] $V_k\subset G_{f,b}$ and
$V_k\in\VB(\delta,f)$ for all $k\in\KC$; \item[(3)]
$\aleph\{P_j\}\leq3C_{d,2}+1$ and $\,\aleph\{V_k\}\leq C_{d,2}$;
\item[(4)] $G_{f,b}\subset\bigcup_{j\in\JC,\,k\in\KC}
\left(\overline{P_j}\bigcup\overline{V_k}\right)$; \item[(5)]
$\#\{k\in\KC\;|\;\mu_d(V_k)\leq 2^{1-d}\,\delta^d\}\leq
C_{d,3}\,\VC_{\delta/2}(f,Q_a^{(d-1)})\,$ and
\newline
$\#\{j\in\JC\;|\;\mu_d(P_j)\leq(2\sqrt{d})^{-d}\,\delta^d\}\leq
C_{d,3}\sum_{m=0}^{m_\delta}2^m\,\VC_{2^{m-1}\delta}(f,Q_a^{(d-1)})\,$,
\newline
where $\,m_\delta:=\min\,\{m\in\Z_+\;|\;2^{m-1}\delta\geq
\Osc(f,Q_a^{(d-1)})\}\,$.
\end{enumerate}
\end{theorem}

\begin{proof}
Let $\{Q^{(d-1)}(k)\}_{k\in\KC_m}$ be the same families of cubes
as in Lemma~\ref{L3.4}, $c_k:=\inf_{x\in Q^{(d-1)}(k)}f(x)$,
$\,b_k=c_k-\delta\,$, $V_k:=G_{f,b_k}(Q^{(d-1)}(k))$ and
$$
P_{m,k,n}:=Q^{(d-1)}(k)\times(c_k-n\delta,c_k-n\delta+\delta)\,,
$$
where $\,k\in\bigcup_m\KC_m\,$ and $n\in\Z_+$. Denote
$\NC_m:=\{2^m+1,\ldots,2^m+2^{m+1}\}\,$. Lemma~\ref{L3.4}(4)
implies that
\begin{equation}
\bigcup_{k\in\KC_m,n\in\NC_m}P_{m,k,n}\ \subset\ \{x\in
G_f\;|\;2^m\delta\leq f(x')-x_d\leq 2^{m+2}\delta\}\,,\label{3.2}
\end{equation}
for all $m=0,1,\ldots,m_\delta$. Let $\;\KC:=\KC_0\,$,
$\,\JC_*:=\bigcup_{m=0}^{m_\delta}\KC_m\times\NC_m\;$ and
$\,\{P_j\}_{j_*\in\JC_*}:=
\bigcup_{m=0}^{m_\delta}\{P_{m,k,n}\}_{k\in\KC_m,n\in\NC_m}\,$.

Assume that $x\in G_f$. If $f(x')-x_d\leq2\delta$ then, by
Lemma~\ref{L3.4}(1), we have
$x\in\bigcup_{k\in\KC}\left(\overline{V_k}\bigcup\overline{P_{0,k,2}}\right)\,$.
If $f(x')-x_d>2^{m_\delta+1}\delta$ then
$$
\dist(x,\Gamma_f)\geq f(x')-x_d-2\,\Osc(f,Q_a^{(d-1)})\geq
f(x')-x_d-2^{m_\delta}\delta>2^{m_\delta}\delta\geq\delta.
$$
Finally, if $2\delta\leq f(x')-x_d\leq 2^{m_\delta+1}\delta$ then
$2^{m+1}\delta\leq f(x')-x_d\leq2^{m+1}\delta+2^m\delta$ for some
nonnegative integer $m\leq m_\delta$ and, in view of
Lemma~\ref{L3.4}(1) and Lemma~\ref{L3.4}(4), we have
$x\in\bigcup_{k\in\KC_m,n\in\NC_m}P_{m,k,n}$. Therefore
\begin{equation}\label{3.3}
\{x\in
G_f\;|\;\dist(x,\Gamma_f)\leq\delta\}\subset\bigcup_{j_*\in\JC_*,\,k\in\KC}
\left(\overline{P_{j_*}}\bigcup\overline{V_k}\right)\,.
\end{equation}

Let us choose a constant
$\,c\in(\delta/(2\sqrt{d}),\delta/\sqrt{d}]\,$ in such a way that
$\,a/c\in\N\,$ and split the set
$\,\overline{Q_a^{(d-1)}}\times[b,+\infty)\,$ into the union of
congruent closed cubes $\,\overline{Q_c^{(d-1)}(i)}\,$ whose
interiors $\,Q_c^{(d-1)}(i)\,$ are mutually disjoint. Let
$\,\{P_j\}_{j\in\JC}\,$ be the collection of all the rectangles
$\,P_{j_*}\,$ and all the cubes $\,Q_c^{(d-1)}(i)\,$ which are
contained in $\,G_{f,b}\,$. Then (1) and (2) are obvious. The
second inequality (3) and (5) follow from the corresponding
statements of Lemma~\ref{L3.4}. The first inequality (3) is a
consequence of (\ref{3.2}), Lemma~\ref{L3.4}(3) and the identity
$\aleph\left\{[2^m,2^{m+2}]\right\}_{i\in\Z_+}=3$. It remains to
prove (4).

Let $\,x\in G_f\,$. If $\,\dist(x,\Gamma_f)\leq\delta\,$ then, by
(\ref{3.3}), either $\,x\in\overline{V_k}\,$ for some
$\,k\in\KC\,$ or $\,x\in \overline{P_{j^*}}\,$ for some
$\,j^*\in\JC^*\,$. Since $\,P_{j_*}\in\PB(\delta)\,$ and
$\,b\leq\inf f-2\delta\,$, in the latter case $\,P_{j_*}\subset
G_{f,b}\,$. If $\,\dist(x,\Gamma_f)>\delta\,$ then the cube
$\,Q_c^{(d-1)}(i)\,$, whose closure contains $\,x\,$, is a subset
of $\,G_{f,b}\,$ because its diameter does not exceed
$\,\delta\,$. Therefore (4) holds true.
\end{proof}

In the two dimensional case we also have the following, more
precise result.

\begin{theorem}\label{T3.6}
Let the conditions of Theorem~{\rm\ref{T3.5}} be fulfilled and
$\,d=2\,$. Then there exists countable families of sets
$\,\{P_j\}_{j\in\JC}\,$ and $\,\{V_k\}_{k\in\KC}\,$ such that
\begin{enumerate}
\item[(1)] $P_j\subset G_{f,b}$ and $P_j\in\PB(\delta)$ for all
$j\in\JC$; \item[(2)] $V_k\subset G_{f,b}$ and
$V_k\in\VB(\delta,f)$ for all $k\in\KC$; \item[(3)]
$\aleph\left(\{P_j\}_{j\in\JC}\bigcup\{V_k\}_{k\in\KC}\right)\leq2$;
\item[(4)] $G_{f,b}\subset\bigcup_{j\in\JC,\,k\in\KC}
\left(\overline{P_j}\bigcup\overline{V_k}\right)$; \item[(5)]
$\#\{k\in\KC\;|\;\mu_2(V_k)\leq\delta^2/2\}\leq
\VC_{\delta/2}(f,Q_a^{(1)})\,$ and
\newline
$\#\{j\in\JC\;|\;\mu_2(P_j)\leq\delta^2/8\}\leq
6\,\VC_{\delta/2}(f,Q_a^{(1)})+12a/\delta\,$.
\end{enumerate}
\end{theorem}

\begin{proof}
In the two dimensional case we do not need Besicovitch's theorem
because the `cube' $\,Q^{(1)}_a\,$ coincides with an interval of
the form $\,(b,b+a)\,$. Given $\,\eps>0\,$, one can easily
construct a finite family $\,\{Q^{(1)}(x)\}_{x\in\XC}\,$ of
disjoint subintervals $\,Q^{(1)}(x)\in(a,a+b)\,$ satisfying the
conditions (1)--(4) of Corollary~\ref{C3.2} with
$\,C_{d,2}=C_{d,3}=1\,$. Therefore Lemma~\ref{L3.4} remains valid
if we substitute $\,C_{d,2}=C_{d,3}=1\,$.

Let $\,k\in\KC:=\XC\,$ and $\,b_k\,$, $\,Q^{(1)}(k)\,$ and
$\,V_k=G_{f,b_k}(Q^{(1)}(k))\,$ be the same as in the proof of
Theorem~\ref{T3.5}. By the above, the first inequality in
Theorem~\ref{T3.5}(5) holds true with $\,C_{d,3}=1\,$. Therefore
$\,\#\KC\leq\VC_{\delta/2}(f,Q_a^{(1)})+2a/\delta\,$ (the second
term is the maximal number of intervals $\,Q^{(1)}(k)\,$ whose
length exceeds $\,\delta/2\,$).

Let $\,V_f:=\bigcup_{k\in\KC}V_k\,$. The set $\,G_f\setminus
V_f\,$ is a polygon with edges parallel to coordinate axes which
has at most $\,2\,\VC_{\delta/2}(f,Q_a^{(1)})\,$ vertices lying on
the horizontal lines $\,\{x\;|\;x_1\in Q_a^{(1)}\,,\,x_2=b_k\}\,$.
Let us choose a constant $\,c\in(\delta/2,\delta]\,$ in such a way
that $\,a/c\in\N\,$ and split the interval $\,Q^{(1)}_a\,$ into
the union of $\,a/c\,$ intervals $\,(a_l,a_{l+1})\,$ of length
$\,c\,$; if $\,a<\delta\,$ then we take
$\,(a_1,a_2):=Q^{(1)}_a\,$. Denote
$$
\KC'_l\ :=\ \{k\in\KC\;|\;[a_{l-2},a_{l+3}]\bigcap
\overline{Q^{(1)}(k)}\ne\emptyset\}\,,\quad
b_{k,\,l}:=\min_{k\in\KC'_l}b_k\,,
$$
and $\,P_{k,\,l}:=(a_l,a_{l+1})\times(b_k,b'_k)\,$ where
$\,b'_k:=\min\{b_{k'}\;|\;b_{k'}>b_k,\,k'\in\KC'_l\}\,$; we assume
that $\,P_{k,\,l}:=\emptyset\,$ whenever
$\,b_k=\max\{b_{k'}\;|\;k'\in\KC'_l\}\,$.

We have $\,\dist(x,\Gamma_f)>\delta\,$ whenever
$\,x_1\in[a_l,a_{l+1}]\,$ and $\,x_2<b_{k,\,l}\,$. Therefore
$$
\{x\in G_f\setminus V_f
\;|\;\dist(x,\Gamma_f)\leq\delta,\,x_1\in[a_l,a_{l+1}]\}\ \subset\
\bigcup_{k\in\KC'_l}\overline{P_{k,\,l}}
$$
and, consequently, (\ref{3.3}) holds true with
$\,\JC_*:=\bigcup_l\KC'_l\,$ and
$\,\{P_{j_*}\}_{j_*\in\JC_*}:=\bigcup_l\{P_{k,\,l}\}_{k\in\KC'_l}\,$.
For each fixed $\,l\,$ the number of rectangles $\,P_{k,\,l}\,$
does not exceed $\,\#\KC'_l-1\,$. We also have
$\,\sum_l(\#\KC'_l-1)\leq6\,(\#\KC)\,$ because each point
$\,x_1\in Q^{(1)}_a\,$ belongs to at most six intervals
$\,[a_{l-2},a_{l+3}]\,$. Therefore
$$
\#\JC_*\ \leq\ 6\,(\#\KC)\ \leq\
6\,\VC_{\delta/2}(f,Q_a^{(1)})+12a/\delta\,.
$$
The rest of the proof repeats that of Theorem~\ref{T3.5}.
\end{proof}

\subsection{General domains}\label{S3.3}
We shall need the following elementary lemma.

\begin{lemma}\label{L3.7}
Let $h$ be a real-valued function on $\R_+$ and $0<a\leq b\,$. If
the function $\,th(t)\,$ is nondecreasing then
$$
\sum_{i\in\Z\;|\;a\,\leq\,2^i\,\leq\,b}h(2^i)\ \leq\
2\int_{a}^{2b}t^{-1}\,h(t)\,\dR t\,.
$$
\end{lemma}

\begin{proof}
We have $\,\sum_{a\,\leq\,2^i\,\leq\,b}h(2^i)
=2\sum_{a\,\leq\,2^i\,\leq\,b}(2^{-i}-2^{-i-1})\,(2^i)\,h(2^i)\,$.
Since the function $\tilde h(s)=s^{-1}\,h(s^{-1})\,$ is
decreasing, the right hand side is estimated by
$\,2\int_{(2b)^{-1}}^{a^{-1}}s^{-1}\,h(s^{-1})\,\dR s
=2\int_{a}^{2b}t^{-1}\,h(t)\,\dR t\,$.
\end{proof}

\begin{corollary}\label{C3.8}
Let $\Omega\in BV_{\tau,\infty}$. Then for each
$\delta\in(0,\delta_\Omega]$ there exist families of sets
$\{P_j\}_{j\in\JC}$ and $\{V_k\}_{k\in\KC}$ satisfying the
following conditions:
\begin{enumerate}
\item[(1)] for each $j$ there exists $\,l\in\LC\,$ such that
$P_j\subset\Omega_l$ and $U_l(P_j)\in\PB(\delta)$; \item[(2)] for
each $k$ there exists $\,l\in\LC\,$ such that $V_k\subset\Omega_l$
and $U_l(V_k)\in\VB(\delta)$; \item[(3)] $\aleph\{P_j\}\leq
n_\Omega\,(3C_{d,2}+1)$ and $\aleph\{V_k\}\leq n_\Omega\,C_{d,2}$;
\item[(4)] $\Omega_{\delta_0}^\bR\
\subset\bigcup\limits_{j\in\JC,\,k\in\KC}
\left(\overline{P_j}\bigcup\overline{V_k}\right)\ \subset\
\Omega_{\delta_1}^\bR\,$, \item[(5)] $\#\KC\leq
C_{d,3}\,C_{\Omega,\,\tau}\,\tau(2/\delta)
+n_\Omega\,C_{d,2}\,2^{d-1}\,\delta^{-d}\,\mu_d(\Omega_{\delta_1}^\bR)\;$
and
$$
\#\JC\ \leq\ 4\,C_{d,3}\, C_{\Omega,\,\tau}\,\delta^{-1}
\int_{(2D_\Omega)^{-1}}^{4/\delta}t^{-2}\,\tau(t)\,\dR t\;
+\;n_\Omega\,(3C_{d,2}+1)\,
(2\sqrt{d})^d\,\delta^{-d}\,\mu_d(\Omega_{\delta_1}^\bR)\,,
$$
\end{enumerate}
where $\delta_0:=\delta/\sqrt{d}$ and
$\delta_1:=\sqrt{d}\,\delta+\delta/\sqrt{d}\,$.
\end{corollary}

\begin{proof}
Let $\,\Omega_l=U_l^{-1}(G_{f_l,\,b_l})\,$ be the sets introduced
in Subsection~\ref{S1.1}. Given $\,\delta\in(0,\delta_\Omega]\,$,
we apply Theorem~\ref{T3.5} for each $l\in\LC$ and denote by
$\{P_j\}_{j\in\JC(l)}$ and $\{V_k\}_{k\in\KC(l)}$ the families of
subsets of $\,\Omega_l\,$, which satisfy the conditions of
Theorem~\ref{T3.5} in an appropriate orthogonal coordinate system.

Let $\,\JC'(l):=\{j\in\JC(l)
\;|\;\dist(P_j,\partial\Omega)\leq\delta_0\}\,$,
$$
\{P_j\}_{j\in\JC}:=\bigcup_{l\in\LC}\{P_j\}_{j\in\JC'(l)}
\quad\text{and}\quad\{V_k\}_{k\in\KC}
:=\bigcup_{l\in\LC}\{V_k\}_{k\in\KC(l)}\,.
$$
Then each of the conditions (1)--(3) is a consequence of the
corresponding condition in Theorem~\ref{T3.5}.

If $\,x\not\in\bigcup_{l\in\LC}\Omega_l\,$ then
$\,\dist(x,\partial\Omega)\geq\delta_\Omega>\delta_0\,$. If
$\,x\in\Omega_l\bigcap\Omega_{\delta_0}^\bR\,$ then, by
Theorem~\ref{T3.5}(4), we have
$\,x\in\bigcup_{j\in\JC(l),\,k\in\KC(l)}\left(\overline{P_j}
\bigcup\overline{V_k}\right)$. In this case
$\,x\in\bigcup_{j\in\JC'(l),\,k\in\KC(l)}\left(\overline{P_j}
\bigcup\overline{V_k}\right)$ because $\,\diam
P_j\leq\sqrt{d}\,\delta\,$. Therefore $\Omega_{\delta_0}^\bR\,$ is
a subset of $\,\bigcup_{j\in\JC,\,k\in\KC}
\left(\overline{P_j}\bigcup\overline{V_k}\right)$. The estimates
$\,\sup_{x\in V_k}\dist(x,\partial\Omega)\leq\sqrt{d}\,\delta\,$
and $\,\diam P_j\leq\sqrt{d}\,\delta\,$  imply the second
inclusion (4).

In order to prove (5), let us denote by $M_\delta$ the smallest
positive integer such that $\,2^{M_\delta-1}\delta\geq
D_\Omega\,$. By Theorem~\ref{T3.5}(5), we have
$$
\#\{j\in\bigcup_{l\in\LC}\JC(l)\;|\;\mu_d(P_j)\leq
2^{1-d}\,\delta^d\}\leq
C_{d,3}\,C_{\Omega,\,\tau}\sum_{m=0}^{M_\delta}2^m\,\tau({(2^{m-1}\delta})^{-1})\,.
$$
Since $2^{M_\delta-1}\delta\leq2D_\Omega\,$, applying
Lemma~\ref{L3.7} with $\,a=(2D_\Omega)^{-1}\delta\,$, $\,b=2\,$
and $\,h(t)=t^{-1}\,\tau(\delta^{-1}t)\,$, we obtain
$$
\#\{j\in\bigcup_{l\in\LC}\JC(l)\;|\;\mu_d(P_j)\leq
2^{1-d}\,\delta^d\}\ \leq\
4\,C_{d,3}\,C_{\Omega,\,\tau}\,\delta^{-1}
\int_{(2D_\Omega)^{-1}}^{4/\delta}t^{-2}\,\tau(t)\,\dR t\,.
$$
Now the second estimate (5) follows from the first inequality (3)
and the second inclusion (4). Similarly, the first estimate (5) is
a consequence of the second inequality (3), the second inclusion
(4) and the first inequality in Theorem~\ref{T3.5}(5).
\end{proof}

\begin{corollary}\label{C3.9}
Let $\Omega\in BV_{\tau,\infty}$ and $\,\Omega\in\R^2\,$. Then for
each $\delta\in(0,\delta_\Omega]$ there exist families of sets
$\{P_j\}_{j\in\JC}$ and $\{V_k\}_{k\in\KC}$ satisfying the
conditions {\rm (1), (2)} and {\rm(4)} of
Corollary~{\rm\ref{C3.8}} such that
\begin{enumerate}
\item[(3$'$)] $\aleph\left(\{P_j\}_{j\in\JC}
\bigcup\{V_k\}_{k\in\KC}\right)\leq2\,n_\Omega\,$;
\item[(5$'$)]
$\,\#\KC\leq C_{\Omega,\tau}\,\tau(2/\delta)+2\,n_\Omega\,
\delta^{-2}\,\mu_2(\Omega_{\delta_1}^\bR)\,$ and
\newline
$\, \#\JC\ \leq\
6\,C_{\Omega,\tau}\,\tau(2/\delta)\;+\;12\,D_\Omega/\delta\;+\;16\,n_\Omega\,
\delta^{-2}\,\mu_2(\Omega_{\delta_1}^\bR)\,$.
\end{enumerate}
\end{corollary}

\begin{proof}
The corollary is proved in the same way as Corollary~\ref{C3.8},
with the use of Theorem~\ref{T3.6} instead of Theorem~\ref{T3.5}.
\end{proof}

Our proof of Theorem \ref{T1.8} is based on the following simple
lemma.

\begin{lemma}\label{L3.10}
Let $\,\Omega\,$ be an arbitrary domain. Then for every $\delta>0$
there exists a family of sets $\,\{M_k\}_{k\in\KC}\,$ satisfying
the following conditions:
\begin{enumerate}
\item[(1)] $\,M_k\subset\Omega\,$ and $M_k\in\MB(\delta)$ for each
$k\in\KC\,$; \item[(2)] $\aleph\{M_j\}=1\,$; \item[(3)]
$\Omega_{\delta_0}^\bR\ \subset\bigcup\limits_{k\in\KC}
\overline{M_k}\ \subset\ \Omega_{\delta_1}^\bR\,$, where
$\delta_0:=\delta/\sqrt{d}$ and
$\delta_1:=\sqrt{d}\,\delta+\delta/\sqrt{d}\,$.
\end{enumerate}
\end{lemma}

\begin{proof}
Consider an arbitrary cover of $\,\R^d\,$ by closed cubes
$\,\overline{Q_\delta^{(d)}(k)}\,$ with disjoint interiors
$\,Q_\delta^{(d)}(k)\,$ and define
$\,\{M_k\}_{k\in\KC}:=\{\Omega\bigcap
Q_\delta^{(d)}(k)\}_{k\in\KC}\,$, where $\,\KC\,$ the set of
indices $\,k\,$ such that $\,\Omega_{\delta_0}^\bR\bigcap
Q_\delta^{(d)}(k)\ne\emptyset\,$.
\end{proof}

\section{Spectral asymptotics}\label{S4}

\subsection{Estimates of the counting function}\label{S4.1}
In this section we shall always assume that
$\,\delta_0:=\delta/\sqrt{d}\,$,
$\,\delta_1:=\sqrt{d}\,\delta+\delta/\sqrt{d}\,$ and denote
\begin{equation}\label{4.1}
R_\Omega(\lambda,\delta_1)\ :=\
3\,(4\sqrt{d})\,C_{d,1}\int_{\delta_1}^\infty
\left(s^{-1}\lambda^{d-1}+s^{-d}\right)\,\dR(\mu_d(\Omega_s^\bR))\,,
\end{equation}
where
$\,\int\left(s^{-1}\lambda^{d-1}+s^{-d}\right)\,\dR(\mu_d(\Omega_s^\bR))\,$
is understood as a Stieltjes integral.

\begin{theorem}\label{T4.1}
If $\,\Omega\in\R^d\,$ is an arbitrary domain and $\,\delta>0\,$
then
\begin{equation}\label{4.2}
N(\Omega,\lambda)-C_{d,W}\,\mu_d(\Omega)\,\lambda^d\ \geq\
-\,R_\Omega(\lambda,\delta_1)-
C_{d,W}\,\mu_d(\Omega_{4\delta_1}^\bR)\,\lambda^d\,,\quad\forall\lambda>0\,,
\end{equation}
and
\begin{equation}\label{4.3}
N_\DR(\Omega,\lambda)-C_{d,W}\,\mu_d(\Omega)\,\lambda^d\ \leq\
R_\Omega(\lambda,\delta_1)\;+\;((4d)^d+2)\,\delta^{-d}\,\mu_d(\Omega_{4\delta_1}^\bR)
\end{equation}
for all $\,\lambda\leq\delta^{-1}\,$. If $\,\Omega\in
BV_{\tau,\infty}\,$ and $\delta\in(0,\delta_\Omega]$ then
\begin{multline}\label{4.4}
N_\NR(\Omega,\lambda)-C_{d,W}\,\mu_d(\Omega)\,\lambda^d\ \leq\
R_\Omega(\lambda,\delta_1)\;+\;(4d)^d\,\delta^{-d}\,\mu_d(\Omega_{4\delta_1}^\bR)\\
+\;C_{d,6}\,n_\Omega\,\delta^{-d}\,\mu_d(\Omega_{\delta_1}^\bR)\;+\;8\,C_{d,3}\,
C_{\Omega,\,\tau}\,\delta^{-1}
\int_{(2D_\Omega)^{-1}}^{4/\delta}t^{-2}\,\tau(t)\,\dR t
\end{multline}
for all
$\,\lambda\leq\min\{1,C^{1/2}_{d,9}\,n_\Omega^{-1/2}\}\,\delta^{-1}\,$.
\end{theorem}

\begin{proof}
Let $\,Q_{2^{-i}}^{(d)}(i,n)\,$ be the Whitney cubes introduced in
Theorem~\ref{T3.3},
$$
\IC_\delta^-:=\{i\in\IC\;|\;\sqrt{d}\,2^{-i}\leq\delta_0/4\}\,,\quad
\IC_\delta^+:= \{i\in\IC\;|\;\sqrt{d}\,2^{-i}>\delta_1\}\,,
$$
$\,\IC_\delta^0:=\IC\setminus(\IC_\delta^+\bigcup\IC_\delta^-)\,$
and $\,\Omega_\delta^\sigma:=\bigcup_{i\in\IC_\delta^\sigma}
\bigcup_{n\in\NC_i}Q_{2^{-i}}^{(d)}(i,n)$, where $\,\sigma=+\,$,
$\,\sigma=0\,$ or $\,\sigma=-\,$. The set
$\,\Omega_\delta^\sigma\,$ are mutually disjoint and
$\,\overline{\Omega}=\overline{\Omega_\delta^+}\bigcup\overline{\Omega_\delta^0}
\bigcup\overline{\Omega_\delta^-}\,$. By virtue of (\ref{3.1}),
\begin{equation}\label{4.5}
\Omega_\delta^-\subset\Omega_{\delta_0}^\bR\,,\quad
\Omega_\delta^0\subset\Omega_{4\delta_1}^\bR\setminus\Omega_{\delta_0/4}^\bR\,,\quad
\Omega\setminus\Omega_{4\delta_1}^\bR\subset\Omega_\delta^+
\subset\Omega\setminus\Omega_{\delta_1}^\bR\,.
\end{equation}
and
\begin{equation}\label{4.6}
\#\NC_i\ \leq\ 2^{i\,d}\left(\mu_d(\Omega_{4\sqrt{d}\,2^{-i}}^\bR)
-\mu_d(\Omega_{\sqrt{d}\,2^{-i}}^\bR)\right)\,,\quad\forall
i\in\IC\,.
\end{equation}
In view of the second inclusion (\ref{4.5}), we have
\begin{equation}\label{4.7}
\sum_{i\in\IC_\delta^0}\#\NC_i\ \leq\
(4\sqrt{d}\,\delta_0^{-1})^d\,\mu_d(\Omega_{4\delta_1}^\bR)\ =\
(4d)^d\,\delta^{-d}\,\mu_d(\Omega_{4\delta_1}^\bR)\,.
\end{equation}
Since
$\,\aleph\{\,[\sqrt{d}\,2^{-i},4\sqrt{d}\,2^{-i}]\,\}_{i\in\Z}=3\,$
and $\,\Omega_s^\bR=\Omega_{D_\Omega}^\bR\,$ for all $\,s\geq
D_\Omega\,$, the inequalities (\ref{4.6}) imply that
\begin{equation}\label{4.8}
\sum_{i\in\IC_\delta^+}((2^i)^{1-d}\lambda^{d-1}+1)\,\#\NC_i\
\leq\ 3\,(4\sqrt{d})\int_{\delta_1}^\infty
(s^{-1}\lambda^{d-1}+s^{-d})\,\dR(\mu_d(\Omega_s^\bR))
\end{equation}
for all $\,\lambda>0\,$.

By Lemma~\ref{L2.1},
\begin{multline}\label{4.9}
N(\Omega,\lambda) -C_{d,W}\,\mu_d(\Omega)\,\lambda^d\\
\geq\
-\,C_{d,W}\,\mu_d(\Omega\setminus\Omega_\delta^+)\,\lambda^d\,
+\,\left(N_\DR(\Omega_\delta^+,\lambda)
-C_{d,W}\,\mu_d(\Omega_\delta^+)\,\lambda^d\right)\,,
\end{multline}
\begin{multline}\label{4.10}
N_\DR(\Omega,\lambda)
-C_{d,W}\,\mu_d(\Omega)\,\lambda^d\\
\leq\
N_{\NR,\DR}(\Omega\setminus\Omega_\delta^+,\partial\Omega,\lambda)
\,+\,\left(N_\NR(\Omega_\delta^+,\lambda)
-C_{d,W}\,\mu_d(\Omega_\delta^+)\,\lambda^d\right).
\end{multline}
and
\begin{multline}\label{4.11}
N_\NR(\Omega,\lambda)
-C_{d,W}\,\mu_d(\Omega)\,\lambda^d\\
\leq\ N_\NR(\Omega\setminus\Omega_\delta^+,\lambda)
\,+\,\left(N_\NR(\Omega_\delta^+,\lambda)
-C_{d,W}\,\mu_d(\Omega_\delta^+)\,\lambda^d\right)
\end{multline}
Lemma~\ref{L2.1} implies that
\begin{multline*}
\sum_{n\in\NC_i,\,i\in\IC_\delta^+}
\left(N_\DR(Q_{2^{-i}}^{(d)}(i,n),\lambda)-C_{d,W}\,
(2^{-i}\lambda)^d\right)\ \leq\
N(\Omega_\delta^+,\lambda)-C_{d,W}\,
\mu_d(\Omega_\delta^+)\,\lambda^d\\
\leq\sum_{n\in\NC_i,\,i\in\IC_\delta^+}
\left(N_\NR(Q_{2^{-i}}^{(d)}(i,n),\lambda)-C_{d,W}\,
(2^{-i}\lambda)^d\right).
\end{multline*}
In view of Lemma~\ref{L2.8}, the right and left hand sides are
estimated from below and above by $\,\pm\,C_{d,1}
\sum_{i\in\IC_\delta^+}\left((2^i)^{1-d}\lambda^{d-1}+1\right)\#\NC_i$.
Therefore, by (\ref{4.8}),
\begin{equation}\label{4.12}
|\,N(\Omega_\delta^+,\lambda)
-C_{d,W}\,\mu_d(\Omega_\delta^+)\,\lambda^d\,|\ \leq\
R_\Omega(\lambda,\delta_1)\,,\qquad\forall\lambda>0\,.
\end{equation}
Since
$\,\Omega\setminus\Omega_{4\delta_1}^\bR\subset\Omega_\delta^+\,$,
the lower bound (\ref{4.2}) is an immediate consequence of
(\ref{4.9}) and (\ref{4.12}).

Assume that $\,\lambda\leq\delta^{-1}\,$. Let
$\,\{M_k\}_{k\in\KC}\,$ be the family of sets introduced in
Lemma~\ref{L3.10} and
$$
\{S_m\}_{m\in\MC_\DR}\ :=\
\{Q_{2^{-i}}^{(d)}(i,n)\}_{n\in\NC_j,\,i\in\IC_\delta^0}
\bigcup\{M_k\}_{k\in\KC}\,.
$$
Lemma~\ref{L3.10}(3) and (\ref{4.5}) imply that
$\,\bigcup_{m\in\MC_\DR}S_m=\Omega\setminus\Omega_\delta^+\,$. In
view of Lemma~\ref{L3.10}(2), we have
$\,\aleph\{S_m\}_{m\in\MC_\DR}\leq2\,$. Consequently, by
Lemma~\ref{L2.2},
$$
N_{\NR,\DR}(\Omega\setminus\Omega_\delta^+,\partial\Omega,\lambda)\
\leq\
\sum_{m\in\MC_\DR}N_{\NR,\DR}(S_m,\Upsilon_m,\sqrt{2}\,\lambda)\,,
$$
where $\,\Upsilon_m=\partial S_m\bigcap\partial\Omega\,$. Since
each set $\,S_m\,$ belongs either to $\,\PB(d^{-1/2}\delta_1)\,$
or to $\,\MB(\delta)\,$, Lemma~\ref{L2.6} implies that
$\,N_\NR(S_m,\Upsilon_m,\sqrt{2}\,\lambda)\leq1\,$. Moreover, if
$\,S_m\in\MB(\delta)\,$ then, in view of Lemma~\ref{L2.6}(3),
$\,N_\NR(S_m,\Upsilon_m,\sqrt{2}\lambda)>0\,$ only if
$\,\mu_d(S_m)\geq\delta^d-4\pi^{-2}\,\delta^{d+2}\lambda^2\,$. By
Lemma~\ref{L3.10}(3), the number of set
$\,M\in\{M_k\}_{k\in\KC}\,$ satisfying this estimate does not
exceed
$$
\left(1-4\pi^{-2}\,\delta^2\lambda^2\right)^{-1}
\delta^{-d}\,\mu_d(\Omega_{\delta_1}^\bR) \ \leq\
2\,\delta^{-d}\,\mu_d(\Omega_{\delta_1}^\bR)
$$
Taking into account (\ref{4.7}), we obtain
$$
N_{\NR,\DR}(\Omega\setminus\Omega_\delta^+,\partial\Omega,\lambda)\
\leq\ (4d)^d\,\delta^{-d}\,\mu_d(\Omega_{4\delta_1}^\bR)+
2\,\delta^{-d}\,\mu_d(\Omega_{\delta_1}^\bR)\,.
$$
This estimate, (\ref{4.10}) and (\ref{4.12}) imply (\ref{4.3}).

In order to prove (\ref{4.4}), let us consider the family of sets
$\{P_j\}_{j\in\JC}$ and $\{V_k\}_{k\in\KC}$ constructed in
Corollary~\ref{C3.8} and define
$$
\{S_m\}_{m\in\MC_\NR}\ :=\
\{Q_{2^{-i}\delta}^{(d)}(i,n)\}_{n\in\NC_j,\,i\in\IC_\delta^0}
\bigcup\{P_j\}_{j\in\JC}\bigcup\{V_k\}_{k\in\KC}\,.
$$
Corollary~\ref{C3.8}(4) and (\ref{4.5}) imply that
$\,\bigcup_{m\in\MC}S_m=\Omega\setminus\Omega_\delta^+\,$. In view
of Corollary~\ref{C3.8}(3), we have $\,\aleph\{S_m\}_{m\in\MC}\leq
n_\Omega\,C_{d,4}^2\,$. Consequently, by Lemma~\ref{L2.2},
$$
N_\NR(\Omega\setminus\Omega_\delta^+,\lambda)\ \leq\
\sum_{m\in\MC_\NR}N_\NR(S_m,n_\Omega^{1/2}\,C_{d,4}\,\lambda)\,.
$$
Since each set $\,S_m\,$ belongs either to $\,\VB(\delta)\,$ or to
$\,\PB(d^{-1/2}\delta_1)\,$, Lemma~\ref{L2.6} implies that
$\,N_\NR(S_m,n_\Omega^{1/2}\,C_{d,4}\,\lambda)=1\,$ whenever
$\,n_\Omega^{1/2}\,C_{d,4}\,\lambda\leq C_{d,5}\,\delta^{-1}\,$.
Estimating $\,\#\MC\,$ with the use of (\ref{4.7}) and
Corollary~\ref{C3.8}(5) and applying the inequalities
$$
(\delta/4)\,\tau(\delta/2)\ =\
\tau(\delta/2)\int_{2/\delta}^{4/\delta}t^{-2}\,\dR t\ \leq\
\int_{2/\delta}^{4/\delta}t^{-2}\,\tau(t)\,\dR t\ \leq\
\int_{(2D_\Omega)^{-1}}^{4/\delta}t^{-2}\,\tau(t)\,\dR t\,,
$$
we see that
\begin{multline}\label{4.13}
N_\NR(\Omega\setminus\Omega_\delta^+,\lambda)\ \leq\ 8\,C_{d,3}\,
C_{\Omega,\,\tau}\,\delta^{-1}
\int_{(2D_\Omega)^{-1}}^{4/\delta}t^{-2}\,\tau(t)\,\dR t\\
+\;(4d/\delta)^d\,\mu_d(\Omega_{4\delta_1}^\bR)\;+\;
C_{d,6}\,n_\Omega\,\delta^{-d}\,\mu_d(\Omega_{\delta_1}^\bR)
\end{multline}
for all $\,\lambda\leq C_{d,7}\,n_\Omega^{-1/2}\,\delta^{-1}\,$.
Now (\ref{4.4}) follows from (\ref{4.11}) and (\ref{4.12}).
\end{proof}

\subsection{Two dimensional domains}\label{S4.2}
If $d=2$, $\,\tau(t)=t\,$ and $\,\delta\asymp\lambda^{-1}\,$ then
the first term on the right hand side of (\ref{4.13}) coincides
with $\,c\,\lambda\,\log\lambda\,$, where $c$ is some constant. On
the other hand, for two dimensional domains with smooth boundaries
we have $\,N_\NR(\Omega_{\lambda^{-1}}^\bR,\lambda)\sim\lambda\,$
as $\lambda\to\infty$ (see, for example, \cite{SV}). The following
lemma gives a refined estimate for
$\,N_\NR(\Omega\setminus\Omega_\delta^+,\lambda)\,$, which does
not contain the logarithmic factor.

\begin{lemma}\label{L4.2}
Let $\,\Omega\subset\R^2\,$, $\,\Omega\in BV_{\tau,\infty}\,$,
$\,\delta\in(0,\delta_\Omega]\,$ and $\,\Omega_\delta^+\,$ be
defined as in Subsection {\rm\ref{S4.1}}. Then for all
$\,\lambda\leq\frac{\sqrt{2}}{3}\,n_\Omega^{-1/2}\delta^{-1}\,$ we
have
\begin{equation}\label{4.14}
N_\NR(\Omega\setminus\Omega_\delta^+,\lambda)\ \leq\
7\,C_{\Omega,\tau}\,\tau(2/\delta)+(64+18\,n_\Omega)\,
\delta^{-2}\,\mu_2(\Omega_{4\delta_1}^\bR)+12\,D_\Omega/\delta\,.
\end{equation}
\end{lemma}

\begin{proof}
Applying the same arguments as in the proof of Theorem~\ref{T4.1}
but using Corollary~\ref{C3.9} instead of Corollary~\ref{C3.8},
one obtains (\ref{4.14}) instead of (\ref{4.13}).
\end{proof}

\subsection{Proof of Theorems \ref{T1.3}, \ref{T1.8}
and Corollary~\ref{C1.5}}\label{S4.3} Integrating by parts in the
Stieltjes integral and changing variables $\,s=t^{-1}\,$, we
obtain
\begin{multline}\label{4.15}
\int_\eps^\infty
(s^{-1}\lambda^{d-1}+s^{-d})\,\dR(\mu_d(\Omega_s^\bR))
\;+\;(\eps^{-1}\lambda^{d-1}+\eps^{-d})\,\mu_d(\Omega_\eps^\bR)\\
= \int_0^{\eps^{-1}}(\lambda^{d-1}+d\,t^{d-1})\,
\mu_d(\Omega_{t^{-1}}^\bR)\,\dR t\,,\qquad\forall\eps>0\,.
\end{multline}
Therefore
$\,\left((4\delta_1)^{-1}\lambda^{d-1}+(4\delta_1)^{-d}\right)
\mu_d(\Omega_{4\delta_1}^\bR)\leq(\lambda^{d-1}+d\,\delta_1^{1-d})
\int_0^{\delta_1^{-1}}\,\mu_d(\Omega_{t^{-1}}^\bR)\,\dR t\,$ and
$\,\int_{\delta_1}^\infty
(s^{-1}\lambda^{d-1}+s^{-d})\,\dR(\mu_d(\Omega_s^\bR))\leq
(\lambda^{d-1}+d\,\delta_1^{1-d})
\int_0^{\delta_1^{-1}}\,\mu_d(\Omega_{t^{-1}}^\bR)\,\dR t\,$.
Applying these inequalities and the estimates
(\ref{4.2})--(\ref{4.4}) with $\,\delta_1^{-1}=\lambda\,$ or
$\,\delta^{-1}=C_{d,8}\,n_\Omega^{1/2}\,\lambda\,$, we obtain
(\ref{1.1}) and (\ref{1.6}). The estimate (\ref{1.2}) is proved in
the same manner, using (\ref{4.14}) instead of (\ref{4.13}).
Finally, since $\,\int_a^bt^{-2}\,\tau(t)\,\dR t\leq
b^{d-2}\int_a^bt^{-d}\,\tau(t)\,\dR t\,$, (\ref{1.3}) is a
consequence of (\ref{1.1}) and the following lemma.

\begin{lemma}\label{L4.3}
If $\Omega\in BV_{\tau,\infty}$ then
$$
\mu_d(\Omega_\eps^\bR)\ \leq
C_{d,2}\,3^d\,n_\Omega\,D_\Omega^{d-1}\eps+
C_{d,3}\,3^d\,C_{\Omega,\,\tau}\,\eps^d\,\tau(\eps^{-1})\,,
\qquad\forall\eps>0\,.
$$
\end{lemma}

\begin{proof}
Assume first that $f$ is a continuous function on the closed cube
$\overline{Q_a^{(d-1)}}$. Let $\{Q^{(d-1)}(x)\}_{x\in\XC}$ be the
same family of cubes as in Corollary~\ref{C3.2},
$\Gamma_f(x):=\{z\in\Gamma_f\;|\;z'\in Q^{(d-1)}(x)\}$ and
$\XC_\eps:=\{x\in\XC\,|\,Q^{(d-1)}(x)\in\PB(\eps)\}$.

If $\,\dist(y,\Gamma_f)\leq\eps\,$ then
$\dist(y,\Gamma_f(x))\leq\eps\,$ for some $x\in\XC$. Therefore
\begin{multline*}
\mu_d\left(\{y\in
Q_a^{(d-1)}\;|\;\dist(y,\Gamma_f)\leq\eps\}\right)\\
\leq\ \sum_{x\in\XC}\mu_d\left(\{y\in
Q_a^{(d-1)}\;|\;\dist(y,\Gamma_f(x))\leq\eps\}\right).
\end{multline*}
The set $\,\{y\in
Q_a^{(d-1)}\;|\;\dist(y,\Gamma_f(x))\leq\eps\}\,$ lies in the
$\,\eps$-neighbourhood of the rectangle
$\,Q^{(d-1)}(x)\times\left(\inf_{z\in
Q^{(d-1)}(x)}f(z)\,,\,\sup_{z\in Q^{(d-1)}(x)}f(z)\right)\,$. In
view of Corollary~\ref{C3.2}(4), the measure of this
$\,\eps$-neighbourhood does not exceed
$\,3\eps\,(a_x+2\eps)^{d-1}\,$, where $a_x$ is the length of the
edge of $Q^{(d-1)}(x)$. Therefore
$$
\mu_d\left(\{y\in
Q_a^{(d-1)}\;|\;\dist(y,\Gamma_f)\leq\eps\}\right)
\leq3^d\,\eps^d\,(\#\XC_\eps)+\sum_{x\in\XC\setminus\XC_\eps}3^d\,
\eps\,a_x^{d-1}\,.
$$
Now  the obvious inequality $ \sum_{x\in\XC}a_x^{d-1}\leq
a^{d-1}\,\aleph\{Q^{(d-1)}(x)\}_{x\in\XC} $ and
Corollary~\ref{C3.2}(3) imply that
$$
\mu_d(\{y\in\R^d\;|\;\dist(y,\Gamma_f)\leq\eps\})\ \leq\
C_{d,2}\,3^d\eps\,a^{d-1}+C_{d,3}\,3^d\eps^d\,\VC_\eps(f,Q_a^{(d-1)})\,.
$$
Since $\,\Omega_\eps^\bR=\bigcup_{l\in\LC}
\{x\in\Omega\,|\,\dist(x,\Gamma_{f_l})\leq\eps\}$, where $\,f_l\,$
are the functions introduced in Subsection~\ref{S1.1}, the lemma
follows from this inequality.
\end{proof}

\subsection{Proof of Corollaries~\ref{C1.6} and \ref{C1.9}}\label{S4.4}
Let $\,\Omega\in\Lip_\alpha\,$, $\,f_l\,$ be the functions
introduced in Subsection {\rm\ref{S1.1}} and
$\,|\Omega|_\alpha:=\max_l|f_l|_\alpha\,$, where
$\,|\cdot|_\alpha\,$ is the seminorm defined in Subsection
\ref{S1.1}. If $\,x\in G_{f_l}\,$ and
$\,\dist(x,(y',f_l(y'))\leq\delta\,$ then
\begin{equation}\label{4.16}
f_l(x')-x_d\ \leq\ |x_d-f_l(y')|+|f_l(y')-f_l(x')|\ \leq\
\delta+\delta^\alpha\,|f_l|_\alpha\,.
\end{equation}
Therefore $\,\{x\in
G_{f_l}\;|\;\dist(x,\Gamma_{f_l})\leq\delta\}\subset\{x\in
G_{g_l}\;|\;
f_l(x')-x_d\leq\delta+\delta^\alpha\,|f_l|_\alpha\}\,$ and,
consequently $\,\mu_d(\{x\in
G_{f_l}\;|\;\dist(x,\Gamma_{f_l})\leq\delta\})\leq
a^{d-1}\,(\delta+\delta^\alpha\,|f_l|_\alpha)\,$. This immediately
implies the following lemma.

\begin{lemma}\label{L4.4}
If $\,\Omega\in\Lip_\alpha\,$ and $\,\delta\leq\delta_\Omega\,$
then $\,\mu_d(\Omega_\delta^\bR)\leq
n_\Omega\,D_\Omega^{d-1}\,(\delta+\delta^\alpha\,|\Omega|_\alpha)\,$.
\end{lemma}

If $\,Q_c^{(d-1)}\subset Q_{a_l}^{(d-1)}\,$ then $\,\diam
Q_c^{(d-1)}=d^{1/2}\,c\,$ and
\begin{equation}\label{4.17}
2\,\Osc(f,Q_c^{(d-1)})\ \leq\ \sup_{x',y'\in
Q_c^{(d-1)}}|f_l(x')-f_l(y')|\ \leq\
d^{\alpha/2}\,c^\alpha\,|f|_\alpha\,.
\end{equation}
Therefore $\,c^{d-1}\geq
d^{(1-d)/2}\,|f|_\alpha^{(1-d)/\alpha}\,\delta^{(d-1)/\alpha}\,$
whenever $\,\Osc(f,Q_c^{(d-1)})\geq\delta/2\,$ and, consequently,
\begin{equation}\label{4.18}
\VC_{\delta/2}(f,Q_a^{(d-1)})\ \leq\
d^{(d-1)/2}\,a^{d-1}\,|f|_\alpha^{(d-1)/\alpha}\,\delta^{(1-d)/\alpha}+1\,.
\end{equation}
The inequality (\ref{4.18}) implies the following result.

\begin{lemma}\label{L4.5}
If $\,\Omega\in\Lip_\alpha\,$ and
\begin{equation}\label{4.19}
\tau(t)\ =\ 2^{(1-d)/\alpha}\,d^{(d-1)/2}\,D_\Omega^{d-1}\,
|\Omega|_\alpha^{(d-1)/\alpha}\,t^{(d-1)/\alpha}+1
\end{equation}
then $\,\Omega\in BV_{\infty,\tau}\,$ and $\,C_{\Omega,\tau}\leq
n_\Omega \,$.
\end{lemma}

Clearly, (\ref{1.4}) follows from (\ref{1.1}) and
Lemma~\ref{L4.5}. Similarly, (\ref{1.6}) and Lemma~\ref{L4.4}
imply (\ref{1.7}). It remains to prove (\ref{1.5}) and
(\ref{1.8}).

Assume that $\,\Omega\in\lip_\alpha\,$. Then for each $\,\eps>0\,$
we can find functions $\,f_{l,\,1}^{(\eps)}\in\Lip_1\,$ and
$\,f_{l,\,2}^{(\eps)}\in\Lip_\alpha\,$ such that
$\,f_l=f_{l,\,1}^{(\eps)}+f_{l,\,2}^{(\eps)}\,$ and
$\,|f_{l,\,2}^{(\eps)}|_\alpha\leq\eps\,$. Obviously,
$\,\VC_\delta(f_{l,\,1}^{(\eps)}+f_{l,\,2}^{(\eps)},Q)\ \leq\
\VC_{\delta/2}(f_{l,\,1}^{(\eps)},Q)\;+\;
\VC_{\delta/2}(f_{l,\,2}^{(\eps)},Q)\,$. Therefore (\ref{4.18})
implies that
\begin{multline}\label{4.20}
\VC_\delta(f_l,Q_{a_l}^{(d-1)})\ \leq\
d^{(d-1)/2}\,D_\Omega^{d-1}\left(\eps^{(d-1)/\alpha}\,\delta^{(1-d)/\alpha}
+C_\eps^{d-1}\,\delta^{1-d}\right)+2\\
\leq\ \eps^{(d-1)/\alpha}\,\tau_\eps(\delta^{-1})\,,
\end{multline}
where $\,C_\eps:=\max_l\,|f_{l,\,2}^{(\eps)}|_1\,$, and
$$
\tau_\eps(t)\ :=\
d^{(d-1)/2}\,D_\Omega^{d-1}\left(t^{(d-1)/\alpha}
+C_{\eps,\Omega}\,\eps^{(1-d)/\alpha}\,t^{d-1}\right)+2\,\eps^{(1-d)/\alpha}\,.
$$
We also have
$$
|f_l(x')-f_l(y')|\leq\eps\,|x'-y'|^\alpha+|f_{l,\,2}^{(\eps)}|_1\,|x'-y'|\,,
\qquad\forall x',y'\in Q_{a_l}^{(d-1)}\,.
$$
Therefore, instead of (\ref{4.16}), we obtain $\,f_l(x')-x_d\leq
\delta+\delta\,|f_{l,\,2}^{(\eps)}|_1+\delta^\alpha\,\eps\,$. This
implies that $\,\mu_d(\Omega_\delta^\bR)\leq
n_\Omega\,D_\Omega^{d-1}\,(\delta+C_{\eps}\,\delta+\delta^\alpha\,\eps)\,$
whenever $\,\delta\leq\delta_\Omega\,$.

In view of (\ref{4.20}), we have $\,\Omega\in
BV_{\infty,\tau_\eps}\,$ and
$\,C_{\Omega,\tau_\eps}\leq\eps^{(d-1)/\alpha}\,n_\Omega\,$.
Choosing a sufficiently large constant $\,C\,$ and applying
(\ref{4.2})--(\ref{4.4}) with $\,\delta=C\,\lambda^{-1}\,$ and
$\,\tau=\tau_\eps\,$, we see that
\begin{align*}
|\,N_\NR(\Omega,\lambda)-C_{d,W}\,\mu_d(\Omega)\,\lambda^d\,|\
&\leq\ \eps^{(d-1)/\alpha}\,C'_\Omega\,\lambda^{(d-1)/\alpha}
+C'_{\Omega,\eps}\,\lambda^{d-1}\,,\\
|\,N_\DR(\Omega,\lambda)-C_{d,W}\,\mu_d(\Omega)\,\lambda^d\,|\
&\leq\ \eps\,C'_\Omega\,\lambda^{d-\alpha}
+C'_{\Omega,\eps}\,\lambda^{d-1}
\end{align*}
for all $\,\lambda>1\,$, where $\,C'_\Omega\,$ is a constant
depending only on the domain $\,\Omega\,$ and
$\,C'_{\Omega,\eps}\,$ is a constant depending on $\,\Omega\,$ and
$\,\eps\,$. Since $\,\eps\,$ can be made arbitrarily small, these
inequalities imply (\ref{1.5}) and (\ref{1.8}). \qed

\subsection{Proof of Theorem~\ref{T1.10}}\label{S4.5}
Let $\,Q_1^{(d-1)}=(0,1)^{d-1}\,$, $\,\alpha\in(0,1)\,$ and
$\,p\,$ be a sufficiently large positive integer. In particular,
we shall be assuming that
$\,p\geq\max\{\alpha^{-1},(1-\alpha)^{-1}\}\,$ and, consequently,
\begin{equation}\label{4.21}
2^{1-\alpha p}\leq1\,,\quad\left(1-2^{-\alpha p}\right)^{-1}\
\leq\ 2\,,\quad \left(1-2^{(1-\alpha)\,p}\right)^{-1}\ \leq\ 2
\end{equation}
and
\begin{equation}\label{4.22}
\left(2^{(1-\alpha)\,(n+1)p}-1\right)\left(2^{(1-\alpha)\,
p}-1\right)^{-1}\ \leq\ 2^{1+(1-\alpha)\,np}\,,\qquad\forall
n=1,2,\ldots
\end{equation}

Given $\,j\in\Z_+\,$, let us denote by $\,\KC_j\,$ the set of
nonnegative integer vectors
$\,\kB=(k_1,\ldots,k_{d-1})\in\Z_+^{d-1}\,$ such that
$\,\max_i\,k_i\leq2^{jp}-1\,$ and consider the
$\,(d-1)$-dimensional cubes
$$
Q(j,\kB)\ :=\ \{x'\in\R^{d-1}\;|\;2^{jp}x'-\kB\in Q_1^{(d-1)}
\}\,,\qquad\kB\in\KC_j\,.
$$
with edges of length $\,{2^{-jp}}\,$. For each fixed
$\,j\in\Z_+\,$ and $\,\kB\in\KC_j\,$ the cubes $\,Q(j,\kB)\,$ are
disjoint and $\,\overline{Q_1^{(d-1)}}=
\bigcup_{\kB\in\KC_j}\overline{Q(j,\kB)}\,$.

Let $\,\psi\in\Lip_1\,$ be a nonnegative Lipschitz function on
$\,Q_1^{(d-1)}\,$ vanishing on the boundary $\,\partial
Q_1^{(d-1)}\,$, $\,a_\psi:=\sup\psi\,$ and
$\,b_{\psi,p}:=\sqrt{d}\,2^{3-(1-\alpha)p}\,(|\psi|_1+a_\psi)\,$.
We shall be assuming that $\,p\,$ is large enough so that
$\,a_\psi>b_{\psi,p}\,$. Let us extend $\,\psi\,$ by 0 to the
whole space $\,\R^{d-1}\,$ and define
$$
g_j(x')\ :=\ \sum_{\kB\in\KC_j} \psi(2^{jp}x'-\kB)\,,\qquad
f_n(x')\ :=\ \sum_{j=0}^n\eps_j\,2^{-\alpha\,jp}\,g_j(x')
$$
and
$\,f(x'):=\lim_{n\to\infty}f_n(x')=\sum_{j=0}^\infty\eps_j\,2^{-\alpha
jp}\,g_j(x')\,$, where $\,\{\eps_j\}\,$ is a nonincreasing
sequence such that $\,\eps_j\in[0,1]\,$ and
\begin{equation}\label{4.23}
2^{(1-\alpha)\left([j/2]-j\right)p}\ \leq\ \eps_{[j/2]}\ \leq\
2\,\eps_j\,,\qquad\forall j=1,2,\ldots
\end{equation}
Note that the condition (\ref{4.23}) is fulfilled whenever
$\,\{\eps_j\}\,$ is a sufficiently slowly decreasing sequence.

\begin{lemma}\label{L4.6}
We have
\begin{enumerate}
\item[(1)]
$\,g_j=0\,$ on $\,\partial Q(j,\kB)\,$ for all
$\,\kB\in\KC_n\,$ and $\,j\geq n\,$;
\item[(2)]
$\,0\leq
f(x')-f_n(x')\leq2\,\eps_{n+1}\,2^{-\alpha\,(n+1)p}\,a_\psi
\leq\eps_{n+1}\,2^{-\alpha\,np}\,a_\psi\,$;
\item[(3)]
$\,|f_n|_\beta\leq2^{1+(\beta-\alpha)np}\,(|\psi|_1+a_\psi)\,$ for
all $\,\beta\in[\alpha,1]\,$;
\item[(4)]
$\,f\in\Lip_\alpha\,$ and
$\,|f|_\alpha\leq2\,(|\psi|_1+a_\psi)\,$;
\item[(5)]
$\,f\in\lip_\alpha\,$ whenever $\,\eps_j\to0\,$ as
$\,j\to\infty\,$;
\item[(6)]
$\,2\,\Osc(f_{n-1},Q(n,\kB))\leq\eps_n\,2^{-\alpha\,np}\,b_{\psi,p}\,$
for all $\,\kB\in\KC_n\,$.
\end{enumerate}
\end{lemma}

\begin{proof}
(1) is obvious and (2) immediately follows from (\ref{4.21}). In
order to prove (3), let us fix $\,\beta\in[\alpha,1]\,$, denote
$\,n':=\max\{j\;|\;2^{-jp}\geq|x'-y'|\}\,$,
$\,n'':=\min\{n,n'\}\,$ and estimate
\begin{multline*}
\sum_{j=0}^n\frac{|g_j(x')-g_j(y')|}{2^{\alpha
jp}\,|x'-y'|^\beta}\ =\
\sum_{j=0}^{n''}\frac{|g_j(x')-g_j(y')|}{2^{\alpha
jp}\,|x'-y'|^\beta}\;+\;
\sum_{j=n''+1}^n\frac{|g_j(x')-g_j(y')|}{2^{\alpha jp}\,|x'-y'|^\beta}\\
\leq\ |\psi|_1\sum_{j=0}^{n''}2^{(1-\alpha)jp}\,|x'-y'|^{1-\beta}
\;+\;a_\psi\sum_{j=n''+1}^n2^{-\alpha jp}\,|x'-y'|^{-\beta}\,.
\end{multline*}
In view of (\ref{4.22}), the first term on the right hand side is
estimated by
$\,|\psi|_1\sum_{j=0}^{n''}2^{(1-\alpha)jp+(1-\beta)np}\leq
2^{1+(\beta-\alpha)np}|\psi|_1\,$. If $\,n\leq n'\,$ then the
second term on the right hand side vanishes; if $\,n>n'\,$ then,
by (\ref{4.21}), it does not exceed
$\,2\,a_\psi\,2^{-\alpha(n''+1)p}|x'-y'|^{-\beta}
\leq2\,a_\psi\,2^{(\beta-\alpha)(n''+1)p}\leq2^{1+(\beta-\alpha)np}a_\psi\,$.
Thus,
\begin{equation}\label{4.24}
\sum_{j=0}^n\frac{|g_j(x')-g_j(y')|}{2^{\alpha
jp}\,|x'-y'|^\beta}\ \leq\
2^{1+(\beta-\alpha)np}\,(|\psi|_1+a_\psi)\,.
\end{equation}
This estimate immediately implies (3) and (4). The inclusion (5)
is also a consequence of (\ref{4.24}) because
$\,|f-f_n|_\alpha\leq\eps_{n+1}\,\sup_{x',y'}
\sum_{j=0}^\infty\frac{|g_j(x')-g_j(y')|}{2^{\alpha
jp}\,|x'-y'|^\alpha}\,$.

Finally, in view of (\ref{4.23}) and (\ref{4.24}), we have
\begin{equation}\label{4.25}
(|\psi|_1+a_\psi)^{-1}|f_j|_1\ \leq\
2^{1+(1-\alpha)[j/2]p}\;+\;\eps_{[j/2]}\,2^{1+(1-\alpha)jp}\ \leq\
\eps_j\,2^{3+(1-\alpha)jp}\,.
\end{equation}
Since $\,\diam Q(n,\kB)=\sqrt{d}\,2^{-np}\,$, (\ref{4.25}) with
$\,j=n-1\,$ implies (6).
\end{proof}

Let $\,\Omega:=G_{f,\,0}\,$,
$\,\Omega_{n,\,\kB}:=\{x\in\Omega\;|\;x'\in Q(n,\kB)\,,\,
x_d\in(f_{n-1}(x'),f(x'))\}\,$,
$\,\Upsilon_{n\,,\kB}:=\partial\Omega_{n,\,\kB}\setminus\partial\Omega\,$
and $\,\Omega_{n-1}\,$ be the interior of
$\,\Omega\setminus\left(\bigcup_{\kB\in\KC_n}\Omega_{n,\,\kB}\right)\,$.

Denote $\,a_{n,\,\kB}:=\sup\limits_{x'\in Q(n,\kB)}f_{n-1}(x')\,$
and consider the function
$$
u_{n,\,\kB}(x)\ :=\ \begin{cases}
\sin\left(2^{\alpha\,np}(x_d-a_{n,\,\kB})/\eps_n\right)\,,
&x_d\geq a_{n-1,\,\kB}\,,\\
0\,,&x_d<a_{n-1,\,\kB}\,,
\end{cases}
$$
on $\,\Omega_{n,\,\kB}\,$. We have $\,u_{n,\,\kB}(x)\in
W^{1,2}(\Omega_{n,\,\kB})\,$ and, in view of Lemma~\ref{L4.6}(1),
$\,u_{n,\,\kB}=0\,$ on $\,\Upsilon_{n,\,\kB}\,$. Applying
Lemma~\ref{L4.6}(2) and Lemma~\ref{L4.6}(6), we see that
\begin{multline*}
\int_{\Omega_{n,\,\kB}}|\nabla u_{n,\,\kB}(x)|^2\,\dR x=
\eps_n^{-2}\,2^{2\alpha\,np}\int_{Q(n,\kB)}\int_0^{f(x')-a_{n,\,\kB}}
\cos^2\left(2^{\alpha\,np}\,x_d/\eps_n\right)\dR x_d\,\dR x'\\
\leq\ \eps_n^{-2}\,2^{2\alpha\,np}\int_{Q(n,\kB)}
\int_0^{\eps_n\,2^{-\alpha\,np}(g_n(x')+a_\psi)}
\cos^2\left(2^{\alpha\,np}x_d/\eps_n\right)\dR x_d\,\dR x'\\
=\ \eps_n^{-1}\,2^{\alpha\,np}\,2^{-(d-1)\,np}\int_{Q_1^{(d-1)}}
\int_0^{\psi(x')+a_\psi} \cos^2x_d\,\dR x_d\,\dR x'
\end{multline*}
and
\begin{multline*}
\int_{\Omega_{n,\,\kB}}|u_{n,\,\kB}(x)|^2\,\dR x\ =\
\int_{Q(n,\kB)}\int_0^{f(x')-a_{n,\,\kB}}
\sin^2\left(2^{\alpha\,np}x_d/\eps_n\right)\dR x_d\,\dR x'\\
\geq\ \int_{Q(n,\kB)}
\int_0^{\eps_n2^{-\alpha\,np}(g_n(x')-b_{\psi,p})}
\sin^2\left(2^{\alpha\,np}x_d/\eps_n\right)\dR x_d\,\dR x'\\
=\ \eps_n\,2^{-\alpha\,np}\,2^{-(d-1)\,np}\int_{Q_1^{(d-1)}}
\int_0^{\psi(x')-b_{\psi,p}} \sin^2x_d\,\dR x_d\,\dR x'\,.
\end{multline*}
Therefore $\;\int_{\Omega_{n,\,\kB}}|\nabla u_{n,\,\kB}(x)|^2\,\dR
x\leq c_{\psi,p}^2\,\eps_n^{-2}\,2^{2\alpha
np}\int_{\Omega_{n,\,\kB}}|u_{n,\,\kB}(x)|^2\,\dR x\,$, where
$$
c_{\psi,p}\ :=\
\left(\frac{\int_{Q_1^{(d-1)}}\int_0^{\psi(x')+a_\psi}\cos^2x_d\;\dR
x_d\,\dR
x'}{\int_{Q_1^{(d-1)}}\int_0^{\psi(x')-b_{\psi,p}}\sin^2x_d\;\dR
x_d\,\dR x'}\right)^{1/2}.
$$
This implies that
$\,N_{\NR,\DR}(\Omega_{n,\,\kB},\Upsilon_{n,\,\kB},\lambda)\geq1\,$
whenever $\,\lambda\geq c_{\psi,p}\,\eps_n^{-1}\,2^{\alpha np}\,$.

Assume that $\,\lambda\in\left[c_{\psi,p}\,\eps_n^{-1}\,2^{\alpha
np},\,c_{\psi,p}\,\eps_{n+1}^{-1}\,2^{\alpha(n+1)p}\right)\,$ and,
using Lemma~\ref{L2.1}, estimate
$$
N_\NR(\Omega,\lambda)\ \geq\ N_\DR(\Omega_{n-1},\lambda)\;+\;
\sum_{\kB\in\KC_n}N_{\NR,\DR}(\Omega_{n,\,\kB},\Upsilon_{n,\,\kB},\lambda)\,.
$$
By the above, the second term on the right hand side is not
smaller than $\,\#\KC_n=2^{(d-1)\,np}\geq(c_{\psi,p}\,2^{\alpha
p})^{(1-d)/\alpha}\,\eps_{n+1}^{(d-1)/\alpha}\,\lambda^{(d-1)/\alpha}\,$.
On the other hand, in view of Theorem~\ref{T1.8}, Lemma~\ref{L4.4}
and Lemma~\ref{L4.6}(3) with $\,\beta=\alpha\,$, we have
$$
N_\DR(\Omega_{n-1},\lambda)\ \geq\
C_{d,W}\,\mu_d(\Omega_{n-1})\,\lambda^d\;
-\;C_d\,(|\psi|_1+a_\psi+1)\,\lambda^{d-\alpha}
$$
for all sufficiently large $\,\lambda\,$. Finally, by
Lemma~\ref{L4.6}(2),
$$
\mu_d(\Omega)\,\lambda^d\;-\;\mu_d(\Omega_{n-1})\,\lambda^d\ \leq\
\eps_n\,2^{-\alpha\,(n-1)p}a_\psi\,\lambda^d\ \leq\
a_\psi\,c_{\psi,p}\,(\eps_n/\eps_{n+1})\,2^{2\alpha
p}\,\lambda^{d-1}\,.
$$
Since $\,\eps_n\leq \eps_{[(n+1)/2]}\leq2\,\eps_{n+1}\,$, the
above estimates imply that
\begin{multline}\label{4.26}
N_\NR(\Omega,\lambda)\ \geq\
C_{d,W}\,\mu_d(\Omega_n)\,\lambda^d\;+\;(c_{\psi,p}\,2^{\alpha
p})^{(1-d)/\alpha}\,\eps_{n+1}^{(d-1)/\alpha}\,\lambda^{(d-1)/\alpha}\\
-\;C_d\,(|\psi|_1+a_\psi+1)\,\lambda^{d-\alpha}
\;-\;C_{d,W}\,a_\psi\,c_{\psi,p}\,2^{2\alpha p+1}\,\lambda^{d-1}
\end{multline}
for all $\,\lambda\in\left[c_{\psi,p}\,\eps_n^{-1}\,2^{\alpha
np},\,c_{\psi,p}\,\eps_{n+1}^{-1}\,2^{\alpha(n+1)p}\right)\,$.

By Lemma~\ref{L4.6}(4), $\,\Omega\in\Lip_\alpha\,$ and we have
$\,(d-1)/\alpha>d-\alpha>d-1\,$. Therefore taking
$\,\eps_0=\eps_1=\dots=1\,$, we obtain a domain satisfying the
conditions of Theorem~\ref{T1.10}(1). If $\,\phi\,$ is a
nonnegative function on $(0,+\infty)$ and $\,\phi(\lambda)\to0\,$
as $\,\lambda\to\infty\,$ then we can choose a sequence
$\,\eps_n\,$ converging to zero and satisfying (\ref{4.23}) in
such a way that the function
$\,\phi(\lambda)\,\lambda^{(d-1)/\alpha}\,$ and the last two terms
in (\ref{4.26}) are estimated by $\,(c_{\psi,p}\,2^{\alpha
p})^{(1-d)/\alpha}\,\eps_{n+1}^{(d-1)/\alpha}\,\lambda^{(d-1)/\alpha}\,$
for all $\,\lambda\in\left[c_{\psi,p}\,\eps_n^{-1}\,2^{\alpha
np},\,c_{\psi,p}\,\eps_{n+1}^{-1}\,2^{\alpha(n+1)p}\right)\,$ and
all sufficiently large $\,n\,$. In view of Lemma~\ref{L4.6}(5),
this proves Theorem~\ref{T1.10}(2). \qed

\section{Remarks and generalisations}\label{S5}

\subsection{Poincar\'e inequality}\label{S5.1}
According to the Poincar\'e inequality,
\begin{equation}\label{5.1}
\int_\Omega |u|^2\,\dR x\leq c_\Omega\int _\Omega |\nabla
u|^2\,\dR x\quad\text{whenever}\ u\in W^{1,2}(\Omega)\ \text{and}\
\int_\Omega u\,\dR x=0,
\end{equation}
where $\,c_\Omega\,$ is a positive constant. By Remark~\ref{R2.4},
the Poincar\'e inequality (\ref{5.1}) on a domain $\,\Omega\,$
holds true if and only if the zero eigenvalue of the Neumann
Laplacian is isolated and
$\,c_\Omega\geq\lambda_{1,\NR}^{-2}(\Omega)\,$.

\begin{lemma}\label{L5.1}
Let $\,\Omega\,$ satisfies {\rm(\ref{5.1})} and
$\,\tilde\Omega\subset\R^d\,$. If there exist an invertible map
$F:\Omega\to\tilde\Omega$ and a constant $\,C_F\,$ such that
$\,|F(x)-F(y)|\leq C_F\,|x-y|\,$ for all $\,x,y\in\Omega\,$ and
$\,|F^{-1}(x)-F^{-1}(y)|\leq C_F\,|x-y|\,$ for all
$x,y\in\tilde\Omega$ then $\,\tilde\Omega\,$ also satisfies
{\rm(\ref{5.1})} with a positive constant
$\,c_{\tilde\Omega}=C_d\,C_F^{2d+2}c_\Omega\,$.
\end{lemma}

\begin{proof}
Let $v\in W^{1,2}(\tilde\Omega)$, $u(x):=v(F^{-1}(x))$ and
$c_u:=\int_\Omega u(x)\,\dR x$. Under the conditions of the lemma
the maps $F$ and $F^{-1}$ are differentiable almost everywhere.
Changing variables and estimating the Jacobians, we obtain
$$
\int_{\tilde\Omega}|v(y)-c_u|^2\,\dR y\ \leq\
C_d\,C_F^d\int_\Omega |u(x)-c_u|^2\,\dR x
$$
and
$$
\int_{\tilde\Omega} |\nabla v(y)|^2\,\dR y\ \geq\
C_d\,C_F^{-d-2}\int_\Omega|\nabla u(x)|^2\,\dR x\,.
$$
These two estimates and the Poincar\'e inequality (\ref{5.1})
imply that
$$
\int_{\tilde\Omega}|v(y)|^2\,\dR y\ \leq\
\int_{\tilde\Omega}|v(y)-c_u|^2\,\dR y\ \leq\
C_d\,C_F^{2d+2}c_\Omega\int_{\tilde\Omega} |\nabla v(y)|^2\,\dR y
$$
whenever $\int_{\tilde\Omega}v\,\dR y=0$.
\end{proof}

Lemma~\ref{L5.1} allows one to extend Theorem~\ref{T1.3} to more
general domains.

\begin{theorem}\label{T5.2}
Assume that there exists a finite collection of domains
$\Omega_l\subset\Omega$ such that
\begin{enumerate}
\item[(a)] $\partial\Omega\subset\bigcup_l\overline{\Omega_l}\,$;
\item[(b$'$)] for each $l$ there exist an invertible map
$F_l:\R^d\to\R^d$ satisfying the conditions of
Lemma~{\rm\ref{L5.1}} such that $F_l(\Omega_l)=G_{f_l,\,b_l}\,$,
where $\,f_l\in BV_{\tau,\infty}(Q_{a_l}^{(d-1)})\,$ and
$\,b_l<\inf f_l\,$;
\item[(c)] $\,a_l\leq D_\Omega\,$ and $\,\sup f_l-b_l\leq
D_\Omega\,$ for all $l\in\LC$.
\end{enumerate}
Then {\rm(\ref{1.1})} holds true.
\end{theorem}

\begin{proof}
Let $\,C_{F_l}\,$ be the constant introduced in Lemma~\ref{L5.1}
and $\,C:=\max_lC_{F_l}\,$. Under conditions of the theorem,
Corollary~\ref{C3.8} remains valid if we replace $\,U_l\,$ with
$\,F_l\,$ and take $\,\delta_n:=C^{-1}\,\delta_n\,$. Since
(\ref{5.1}) is equivalent to the identity
$\,N_\NR(\Omega,c_\Omega^{-2})=1\,$, Lemma~\ref{L2.6} and
Lemma~\ref{L5.1} imply that $\,N_\NR(S_m,\lambda)=1\,$ for all
$\,\lambda\leq c'_\Omega\,\delta^{-1}\,$, where $\,S_m\,$ are the
same sets as in the proof of Theorem~\ref{T4.1} and
$\,c'_\Omega\,$ is a constant depending on the domain
$\,\Omega\,$. Therefore, using the same arguments as in
Subsection~\ref{S4.1}, we obtain the estimates (\ref{4.2}) and
(\ref{4.4}) with some other constants (which may depend on
$\,\Omega\,$). In the same way as in Subsection~\ref{S4.3}, these
estimates imply (\ref{1.1}).
\end{proof}

The following example shows that Theorem~\ref{T5.2} is not just a
formal generalization of Theorem~\ref{T1.3}.

\begin{example}\label{E5.3}
Let $\,f\,$ be a nowhere differentiable $\,\Lip_\alpha$-function
on the interval $\,[0,1]\,$. Assume that $\,f>1\,$ and consider
the domain
$$
\Omega\ :=\
\{(\varphi,r)\in\R^2\;|\;\varphi\in(0,1)\,,\,1<r<f(\varphi)\}\,,
$$
where $\,(\varphi,r)\,$ are the polar coordinates on $\,\R^2\,$.
If $\,y_1=r\sin\varphi\,$ and $\,y_2=r\cos\varphi\,$ are the
standard Cartesian coordinates on $\,\R^2\,$ then the map which
takes the point with polar coordinates $\,(\varphi,r)\,$ into the
point with Cartesian coordinates $\,(y_1,y_2)=(\varphi,r)\,$
satisfies the conditions of Lemma~\ref{L5.1}. Therefore, by
Theorem~\ref{T5.2}, we have (\ref{1.1}).

On the other hand, if $\,(x_1,x_2)\,$ are arbitrary Cartesian
coordinates on $\,\R^2\,$ then
$\,x_1(\varphi,r)=r\sin(\varphi+\varphi_0)\,$ and
$\,x_2(\varphi,r)=r\cos(\varphi+\varphi_0)\,$ for some
$\,\varphi_0\in[0,2\pi)\,$. For every subinterval
$\,(a,b)\subset(0,1)\,$ there exist at least two different points
$\,\varphi_1,\varphi_2\in(a,b)\,$ such that
$\,x_1(\varphi_1,f(\varphi_1))=x_1(\varphi_2,f(\varphi_2))\,$
(otherwise the function $\,x_1(\varphi,f(\varphi))\,$ would be
monotone on $\,(a,b)\,$ and, by Lebesgue's theorem, almost
everywhere differentiable). Since
$\,x_2(\varphi_1,f(\varphi_1))\ne x_2(\varphi_2,f(\varphi_2))\,$,
we see that the set $\,\{r=f(\varphi)\}\,$ cannot be represented
as the graph of a continuous function in Cartesian coordinates.

Nowhere differentiable functions $\,f\in\Lip_\alpha\,$ do exist.
For instance, the function
$\,f(t):=\sum_{n=0}^\infty10^{-n}\,\dist(10^nt,\Z)\,$ is not
differentiable at each $\,t\in\R\,$ (see \cite{W} or \cite{RS-N},
Chapter 1, Section 1) but $\,f\in\Lip_\alpha(\R)\,$ for all
$\,\alpha\in(0,1)\,$.
\end{example}

\subsection{Higher order operators}\label{S5.2}
Let us consider, instead of the Laplacian, a homogeneous elliptic
nonnegative operator $\,A(D_x)\,$ of degree $\,2m\,$ with real
constant coefficients and denote by $\,Q_A\,$ its quadratic form
(we use the standard notation $\,D_x:=-i\,\partial_x)\,$. Let
$\,W^{m,2}(\Omega)\,$ be the Sobolev space,
$\,W_0^{m,2}(\Omega)\,$ be the closure of $\,C_0^\infty\,$ in
$\,W^{m,2}(\Omega)\,$ and $\,A_\NR\,$ and $\,A_\DR\,$ be the
self-adjoint operators in the space $\,L^2(\Omega)\,$ generated by
the quadratic form $\,Q_A\,$ with domains $\,W^{m,2}(\Omega)\,$
and $\,W_0^{m,2}(\Omega)\,$ respectively. Then the results of
Section~\ref{S2} remain valid with the following modifications.
\begin{enumerate}
\item[(i)]
In the definitions of $\,N_{\NR,\DR}\,$, $\,N_\NR\,$, $\,N_\DR\,$
and in Lemma~\ref{L2.2} we replace the Dirichlet form
$\,\int_\Omega|\nabla u|^2\,\dR x\,$ with $\,Q_A\,$,
$\,W^{1,2}(\Omega)\,$ with $\,W^{m,2}(\Omega)\,$, $\,\lambda^2\,$
with $\,\lambda^{2m}\,$, and $\,\varkappa^{-1/2}\,$ with
$\,\varkappa^{-1/(2m)}\,$.
\item[(ii)]
The kernel of the operator $\,A_\NR\,$ is the space
$\,\PC_m(\Omega)\,$ of all polynomials on $\,\Omega\,$ whose
degree is strictly smaller than $\,m\,$. Therefore we have
$\,\int_\Omega|u(x)|^2\,\dR x\leq \lambda^{-2m}\,Q_A[u]\,$ for all
$\,u\in W^{1,2}(\Omega)\ominus\PC_m(\Omega)\,$ if and only if
$\,\lambda_{1,\NR}(\Omega)\geq\lambda\,$, where
$\,\lambda_{1,\NR}(\Omega)\,$ is the first nonzero eigenvalue of
$\,A_\NR\,$. If $\,p_u\,$ is the projection of $\,u\in
L_2(\Omega)\,$ onto the subspace $\,\PC_m(\Omega)\,$ then
$\,\|u-p_u\|_{L^2(\Omega)}\leq\|u-p\|_{L^2(\Omega)}\,$ for all
$\,p\in\PC_m(\Omega)\,$ (cf. Remark~\ref{R2.4}).
\item[(iii)]
Let $\,C_{A,W}:=(2\pi)^{-d}\,\mu_d\{\xi\in\R^d:A(\xi)<1\}\,$. Then
there exists a constant $\,C_{A,Q}\,$ such that
$$
-\,C_{A,Q}\,(\delta\lambda)^{d-1}\leq
N(Q_\delta^{(d)},\lambda)-C_{A,W}\,(\delta\lambda)^d\leq
C_{A,Q}\,(\delta\lambda)^{d-1},\quad\forall\lambda>\delta^{-1}\,,
$$
for all $\,\delta>0\,$ (see Remark~\ref{R2.9}).
\item[(iv)]
Instead of Lemma~\ref{L2.6} we have the following result.
\end{enumerate}

\begin{lemma}\label{L5.4}
There exists a constant $\,c_A\,$ depending only on the operator
$A$ and the dimension $d$ such that the following statements hold
true.
\begin{enumerate}
\item[(1)] If $\,P\in\PB(\delta)$ then
$N_\NR(P,\lambda)=\dim\PC_m\,$ for all $\,\lambda\leq
c_A\,\delta^{-1}$. \item[(2)] If $\,V\in\VB(\delta)$ then
$N_\NR(V,\lambda)=\dim\PC_m\,$ for all $\,\lambda\leq
c_A\,\delta^{-1}$. \item[(3)] If $\,M\in\MB(\delta)\,$,
$\,M\subset Q_\delta^{(d)}\,$ and $\,\Upsilon:=\partial M\bigcap
Q_\delta^{(d)}\,$ then we have
$\,N_{\NR,\DR}(M,\Upsilon,\lambda)\leq\dim\PC_m\,$ for all
$\,\lambda\leq c_A\,\delta^{-1}\,$ and
$\,N_{\NR,\DR}(M,\Upsilon,\lambda)=0\,$ for all $\,\lambda\leq
(1-c_A^{-1}\,\delta^{-d}\mu_d(M))_+^{1/(2m)}\,c_A\,\delta^{-1}\,$.
\end{enumerate}
\end{lemma}

\begin{proof}
We shall denote by $\,C\,$ various constants depending only on $A$
and $d$.

Since $\,A(\xi)\leq C\,\sum_{j=1}^d\,\xi_j^{2m}\,$, it is
sufficient to prove the lemma assuming that
$\,A(D_x)=A_m(D_x):=\sum_{j=1}^d\,D_{x_j}^{2m}\,$. Then (1) is
easily obtained by separation of variables. If $\,u\in
W^{m,2}(Q_\delta^{(d)})\,$, $\,u\equiv0\,$ outside $\,M\,$ and
$\,p_u\,$ is the projection of $\,u\,$ onto the subspace
$\,\PC_m(M)\,$ then
\begin{multline*}
\int_M|p_u|^2\,\dR x\ \leq\ \mu_d(M)\sup_{x\in
Q_\delta^{(d)}}|p_u(x)|^2\ \leq\
C\,\mu_d(M)\,\delta^{-d}\int_{Q_\delta^{(d)}}|p_u|^2\,\dR x\\
=\ C\,\mu_d(M)\,\delta^{-d}\left(\int_M|p_u|^2\,\dR
x+\int_{Q_\delta^{(d)}}|u-p_u|^2\,\dR x\right)\,.
\end{multline*}
Applying (ii) and this estimate instead of Remark~\ref{R2.4} and
(\ref{2.6}), we obtain (3) in the same way as Lemma~\ref{L2.6}(3).

In order to prove (2), let us assume that
$\,V=G_{f,\,b}(Q_c^{(d-1)})\,$ with $\,c\leq\delta\,$, $\,b=\inf
f-\delta$ and $\,\Osc f\leq\delta/2\,$ and consider a function
$\,u\in W^{m,2}(V)\,$. Let $\,p_{u;\,r,\,k}(x')\,$ be the
projection of the function $\,\partial_{x_d}^k\,u(x',r)\in
L^2(Q_c^{(d-1)})\,$ onto the subspace
$\,\PC_{m-k}(Q_c^{(d-1)})\,$, $\,p_{u;\,r}(x):=
\sum_{k=0}^{m-1}\,\frac1{k!}\,(x_d-r)^k\;p_{u;\,r,\,k}(x')\,$ and
$\,v_r(x):=\sum_{k=0}^{m-1}\,\frac1{k!}\,(x_d-r)^k\,\partial_{x_d}^ku(x',r)\,$,
where $\,r\in[b,b+\delta]\,$ and $\,x_d\in[b,f(x')]\,$. We have
\begin{equation}\label{5.2}
|u(x)-p_{u;\,r}(x)|^2\ \leq\
2\,|u(x)-v_r(x)|^2+2\,|v_r(x)-p_{u;\,r}(x)|^2\,.
\end{equation}
Since $\,|x_d-b|\leq2\delta\,$, Jensen's inequality implies that
\begin{multline*}
|u(x)-v_r(x)|^2\ =\
((m-1)!)^{-2}\,|\int_r^{x_d}(x_d-t)^{m-1}\,\partial_{x_d}^mu(x',t)\,\dR
t\,|^2\\
\leq\
((m-1)!)^{-2}\,|x_d-r|\,\int_r^{x_d}(x_d-t)^{2m-2}\,|\partial_{x_d}^mu(x',t)|^2\,\dR
t\\
\leq\
((m-1)!)^{-2}\,(2\delta)^{2m-1}\int_b^{f(x')}|\partial_{x_d}^mu(x)|^2\,\dR
x_d\,.
\end{multline*}
In view of (ii) and (1), we also have
$$
\int_{Q_c^{(d-1)}}|\partial_{x_d}^k\,u(x)-p_{u;\,r,\,k}(x')|^2\,\dR
x'\ \leq\ C\,\delta^{2m-2k}\,Q_{A'_{m-k}}[\partial_{x_d}^k\,u(x)]
$$
for all $\,k=0,\ldots,m-1\,$, where
$\,A'_{m-k}(D_{x'}):=\sum_{j=1}^{d-1}\,D_{x_j}^{2m-2k}\,$ and
$\,Q_{A'_{m-k}}\,$ is the quadratic form of $\,A'_{m-k}\,$ with
domain $\,W^{m-k,\,2}(Q_c^{(d-1)})\,$. Therefore, integrating
(\ref{5.2}) over $\,r\in[b,b+\delta]\,$, $\,x_d\in[b,f(x')]\,$,
$\,x'\in Q_c^{(d-1)}$ and estimating $\,|x_d-r|\leq2\delta\,$, we
obtain
\begin{multline*}
\delta^{-1}\int_b^{b+\delta}\int|u(x)-p_{u;\,r}(x)|^2\,\dR
x\,\dR r\\
\leq\ C\,\delta^{2m}\int_V|\partial_{x_d}^mu(x)|^2\,\dR x\
+C\,\delta^{2m}\sum_{k=0}^{m-1}\,\sum_{|\alpha|=m}\int_P
|\partial_x^\alpha u(x)|^2\,\dR x\,,
\end{multline*}
where $\,P=Q_c^{(d-1)}\times(b,b+\delta)\,$. Since the
$\,L_2$-norms of the mixed derivatives $\,\partial_x^\alpha
u(x)\,$ on a rectangle are estimated by the $\,L_2$-norms of the
derivatives $\,\partial_{x_j}^m\,$, this estimate and (ii) imply
(2).
\end{proof}

Applying the same arguments as in Section~\ref{S4} and using (iii)
and Lemma~\ref{L5.4}, we obtain the following result.

\begin{theorem}\label{T5.5}
Let $\,A\,$ be a homogeneous nonnegative elliptic differential
operator of order $\,2m\,$ with real constant coefficients. If
$\,N_\NR(\lambda,\Omega)\,$ and $\,N_\DR(\lambda,\Omega)\,$ denote
the number of eigenvalues of the corresponding self-adjoint
operator lying below $\,\lambda^{2m}\,$ then Theorems
{\rm\ref{T1.3}}, {\rm\ref{T1.8}} and Corollaries {\rm\ref{C1.5}},
{\rm\ref{C1.6}}, {\rm\ref{C1.9}} holds true with
$\,C_{d,W}:=C_{A,W}\,$.

\end{theorem}

\subsection{Other function spaces}\label{S5.3}
Let $\,B^\alpha_{p,q}\,$ be the Besov space and
$\,BV_{\beta,\infty}:=BV_{\tau_\beta,\infty}\,$ where
$\,\tau_\beta(t)=(t^\beta+1)\,$ and $\beta\in(0,+\infty)$.
Lemma~\ref{L4.5} implies that
$\,B^\alpha_{\infty,\infty}=\Lip_\alpha\subset
BV_{(d-1)/\alpha,\infty}\,$. Estimating the norm of the embedding
$\,B_{p,\infty}^\alpha(Q_a^{(d-1)})\hookrightarrow
C(Q_a^{(d-1)})\,$ for $\,\alpha p> d-1\,$ and $\,a>0\,$, one can
also show that $\,B_{p,\infty}^\alpha\subset
BV_{(d-1)/\alpha,\infty}\,$ whenever $\,\alpha p>d-1\,$.

\subsection{Open problems}
\subsubsection{The spaces $BV_{\tau,\infty}$}\label{S5.3.1}
The space $\,BV_{\beta,\infty}\,$ or $\,BV_{\tau,\infty}$ (under
certain conditions on the function $\tau$) is a Banach space with
respect to an appropriate norm. Similar spaces have been
considered in the dimension one, but we could not find references
in the multidimensional case. It would be interesting to find a
more constructive description of these spaces and to investigate
their properties.

\subsubsection{More general domains}\label{S5.3.2}
The crucial point in our proof of Theorem~\ref{T1.3} is the
construction of the families $\,\{S_m\}_\MC\,$ such that
\begin{enumerate}
\item[(i)]
$\,\Omega_\delta^\bR\subset\bigcup_mS_m\subset\Omega\,$,
\item[(ii)]
$\,\aleph\{S_m\}_\MC\leq C\,$,
\item[(iii)]
$\,N_\NR(S_m,\lambda)\leq C'\,$ whenever $\,\lambda\leq
C''\delta^{-1}\,$,
\end{enumerate}
where $\,C\,$, $\,C'\,$ and $\,C''\,$ are some constants
independent of $\,\delta\in\R_+\,$.

The remainder estimate in the Weyl formula for the Neumann
Laplacian depends on the behaviour of $\,\#\MC\,$ as
$\,\delta\to0\,$. In this paper we were assuming that $\,\Omega\,$
is the union of subgraphs of continuous functions, used
Lemma~\ref{L2.6} in order to prove (iii) and applied
Corollary~\ref{C3.2} in order to estimate $\,\aleph\{S_m\}\,$ and
$\,\#\MC\,$. Theorem~\ref{T3.1} allows one to construct families
of open sets $\,S_m\,$ satisfying (i)--(iii) for many other
domains $\,\Omega\,$. It should be possible to find less
restrictive sufficient conditions which guarantee the existence of
such families and imply an asymptotic formula for
$\,N_\NR(\Omega,\lambda)\,$.

\subsubsection{Operators with variable coefficients}\label{S5.3.3}
Our main goal was to estimate the contribution of
$\,\partial\Omega\,$ to the Weyl formula. In the interior part of
$\,\Omega\,$ we used the old fashioned variational technique based
on the Whitney decomposition and Dirichlet--Neumann bracketing.
There are much more advanced methods of studying the asymptotic
behaviour of the spectral function at the interior points (see the
monographs \cite{Iv3}, \cite{SV} or the recent papers \cite{BI},
\cite{Iv4}), which are applicable to operators with variable
coefficients.

Freezing the coefficients at an arbitrary point $\,x\in S_m\,$, we
see that (iii) remains valid for a uniformly elliptic operator
$\,A\,$ with variable coefficients, provided that the
corresponding quadratic form is homogeneous, the coefficients are
uniformly continuous, $\,\delta\,$ is sufficiently small and
$\,\diam S_m\leq c\,\delta\,$ with some constant $\,c\,$
independent of $\,\delta\,$. Using this observation and applying a
more powerful technique in the interior of $\,\Omega\,$, one can
try to extend our results to operators with variable coefficients.

\subsubsection{Reminder estimate for the Dirichlet Laplacian}\label{S5.3.4}
It is not difficult to construct a bounded domain $\,\Omega\,$
such that
$\,\lim_{\delta\to0}\,|\delta^{-\alpha}\,\mu_d(\Omega_\delta^\bR)|=C'\,$
and
\begin{equation}\label{5.3}
N_\DR(\Omega,\lambda) -C_{d,W}\,\mu_d(\Omega)\,\lambda^d\ \geq\
-\,C^{-1}\,\lambda^{d-\alpha}\,,\qquad\forall\lambda>C\,,
\end{equation}
where $C$ and $C'$ are some positive constants. For example, it
can be done by considering a cube with a sequence of `cracks'
converging to the outer boundary, which get denser as the outer
boundary is approached (similar domains were studied in \cite{LV}
and \cite{MV}). For such a domain the estimate (\ref{1.7}) is
order sharp. It would be interesting find a domain
$\,\Omega\in\Lip_\alpha\,$ satisfying (\ref{5.3}) (cf. Theorem
\ref{T1.10}). Note that in the known examples disproving the
so-called Berry conjecture (see, for instance, \cite{BLe} or
\cite{LV}) the domain does not belong to the class
$\,\Lip_\alpha\,$.

\section{Constants}

Throughout the paper $\,C_{d,W}\,$ is the Weyl constant (see
Subsection~\ref{S1.1}),
$$
C_{d,1}:=\sum_{n=0}^{d-1}\frac{n!\,(d-n)!}{d!}\,C_{n,W}\,,\qquad
C_{0,W}:=1\,,
$$
$\,C_{d,2}=2^{d-1}\,\CC_{d-1}\,$ and
$\,C_{d,3}=6^{d-1}\,\hat\CC_{d-1}\,$ where $\,\CC_{d-1}\,$ and
$\,\hat\CC_{d-1}\,$ are the constants introduced in
Theorem~\ref{T3.1},
$$
C_{d,4}:=(4\,C_{d,2}+2)^{1/2}\,,\quad
C_{d,5}:=\min\left\{(1+2\pi^{-2})^{-1/2},\,\pi(1+d^{-1})^{-1}\right\}\,,
$$
$$
C_{d,6}:=2^{d-1}\,C_{d,2}+(3\,C_{d,2}+1)\,(2\sqrt{d})^d\,,\quad
C_{d,7}:=C_{d,4}^{-1}\,C_{d,5}\,,
$$
$$
C_{d,8}:=\max\{1,C_{d,7}^{-1/2}\}\,,\quad
C_{d,9}:=8\,C_{d,3}\,C_{d,8}\,,
$$
$$
C_{d,10}:=(d+1)\left(12\sqrt{d}\,C_{d,1}+4\,C_{d,W}
+(4^dd^d+C_{d,6})\,(4d^{1/2}+4d^{-1/2})^d\right),
$$
$$
C_{d,11}:=(d+1)\left(12\sqrt{d}\,C_{d,1}+4\,C_{d,W}
+(4^dd^d+2)\,(4d^{1/2}+4d^{-1/2})^d\right).
$$

\begin{remark}
If $\rho$ is continuous then  Theorem~\ref{T3.1} holds true with
$\,\CC_n=2^n\,$ and $\,\hat\CC_n=4^n\,$ (see \cite{G}). Since the
function $\rho$  in the proof of Corollary~\ref{C3.2} is
continuous, all our results remain valid for $\,C_{d,2}=4^{d-1}\,$
and $\,C_{d,3}=24^{d-1}$.
\end{remark}


\begin{thebibliography}{MMM}

\bibitem[BD]{BD} V. Burenkov and E.B. Davies. {\it Spectral
stability of the Neumann Laplacian,} J. Differential Equations
{\bf 186} (2002), no. 2, 485--508.

\bibitem[BI]{BI}
M. Bronstein and V. Ivrii. {\it Sharp spectral asymptotics for
operators with irregular coefficients. I. Pushing the limits,}
Comm. Partial Differential Equations {\bf 28} (2003), no. 1-2,
83--102.

\bibitem[BLe]{BLe}
M. van den Berg and M. Levitin. {\it Functions of Weierstrass type
and spectral asymptotics for iterated sets,} Quart. J. Math.
Oxford Ser. {\bf 47} (1996), no. 188, 493--509.

\bibitem[BLi]{BLi}
M. van den Berg and M. Lianantonakis. {\it Asymptotics for the
spectrum of the Dirichlet Laplacian on horn-shaped regions}.
Indiana Univ. Math. J. 50 (2001), no. 1, 299--333.

\bibitem[BS]{BS}
M.S. Birman and M.Z. Solomyak. {\it The principal term of spectral
asymptotics for ``non-smooth'' elliptic problems,} Funktsional.
Anal. i Prilozhen. {\bf 4:4} (1970), 1--13  (Russian), English
transl. in Functional Anal. Appl. {\bf 4} (1971).

\bibitem[F]{F}
 B. Fedosov. {\it Asymptotic formulae for the eigenvalues of the
Laplace operator in the case of a polygonal domain} (Russian),
Dokl. Akad. Nauk SSSR {\bf 151} (1963), 786--789.

\bibitem[G]{G}
M. de Guzm\'an. {\it Differentiation of integrals in $\R^n\,$,}
Springer Verlag, 1975 (Lecture Notes in Mathematics, v. 481).

\bibitem[HSS]{HSS}
R. Hempel, L. Seco and B. Simon. {\it The essential spectrum of
Neumann Laplacians on some bounded singular domains,} J. Funct.
Anal. {\bf 102} (1991), no. 2, 448--483.

\bibitem[Iv1]{Iv1}
V. Ivrii. {\it On the second term of the spectral asymptotics for
the Laplace-Beltrami operator on manifolds with  boundary,}
Funktsional.  Anal. i Prilozhen. {\bf 14} (1980), no. 2, 25--34.
English transl. in Functional Anal. Appl. {\bf 14} (1980),
98--106.

\bibitem[Iv2]{Iv2}
\bysame. {\it The asymptotic Weyl formula for the Laplace-Beltrami
operator in Riemannian polyhedra and domains with conical
singularities of the boundary,} Dokl. Akad. Nauk SSSR 288 (1986),
35--38 (Russian).

\bibitem[Iv3]{Iv3}
\bysame. {\it Microlocal Analysis and Precise Spectral
Asymptotics,} Springer-Verlag, SMM, 1998.

\bibitem[Iv4]{Iv4}
\bysame. {\it Sharp spectral asymptotics for operators with
irregular coefficients. II. Domains with boundaries and
degenerations,} Comm. Partial Differential Equations {\bf 28}
(2003), no. 1-2, 103--128.

\bibitem[LV]{LV}
M. Levitin and D. Vassiliev. {\it Spectral asymptotics, renewal
theorem, and the Berry conjecture for a class of fractals,} Proc.
London Math. Soc. {\bf 72} (1996), 188--214.

\bibitem[M1]{M1}
V. Maz'ya. {\sl On Neumann's problem for domains with irregular
boundaries,} Siberian Math. J. {\bf 9} (1968), 990--1012.

\bibitem[M2]{M2}
V. Maz'ya. {\sl Sobolev spaces,} Leningrad University, Leningrad,
1985. English translation in {\sl Springer Series in Soviet
Mathematics,} Springer-Verlag, Berlin, 1985.

\bibitem[Ma]{Ma}
C. Mason. {\sl Log-Sobolev inequalities and regions with exterior
exponential cusps,} J. Funct. Anal. 198 (2003), 341--360.

\bibitem[Me]{Me}
G. M\'etivier. {\it Valeurs propres de problèmes aux limites
elliptiques irr\'eguli\`eres} (French), Bull. Soc. Math. France
Suppl. M\'em. {\bf 51--52} (1977), 125--219.

\bibitem[Mi]{Mi}
Y. Miyazaki. {\it A sharp asymptotic remainder estimates for the
eigenvalues of operators associated with strongly elliptic
sesquilinear forms,} Japan. J. Math. {\bf 15} (1989), no. 1,
65--97.

\bibitem[MV]{MV}
S. Molchanov and B. Vainberg. {\it On spectral asymptotics for
domains with fractal boundaries of cabbage type,} Math. Phys.
Anal. Geom. {\bf 1} (1998), no. 2, 145--170.

\bibitem[RS-N]{RS-N}
F. Riesz and B. Sz.-Nagy. {\it Le\c{c}ons d'analyse
fonctionnelle,} (French) Académie des Sciences de Hongrie,
Akadémiai Kiadó, Budapest, 1952. English translation: {\it
Functional analysis,} Dover Publications Inc., New York, 1990.

\bibitem[Se]{Se}
R. Seeley. {\it An estimate near the boundary for the spectral
function of the Laplace operator,} Amer. J. Math. {\bf 102}
(1980), no. 5, 869--902.

\bibitem[St]{St}
E. Stein. {\it Singular integrals and differentiability properties
of functions,} Princeton University Press, Princeton, 1970.

\bibitem[Sa]{Sa}
Yu. Safarov. {\it Fourier Tauberian Theorems and applications,} J.
Funct. Anal. {\bf 185} (2001), 111--128.

\bibitem[Si]{Si}
B. Simon. {\it The Neumann Laplacian of a jelly roll,} Proc. Amer.
Math. Soc. {\bf 114} (1992), no. 3, 783--785.

\bibitem[SV]{SV}
Yu. Safarov and D. Vassiliev. {\it The asymptotic distribution of
eigenvalues of partial differential operators,} American
Mathematical Society, 1996.

\bibitem[V]{V}
D. Vassiliev. {\it Two-term asymptotic behavior of the spectrum of
a boundary value problem in the case of a piecewise smooth
boundary,} Dokl. Akad. Nauk SSSR {\bf 286} (1986), no. 5,
1043--1046. English translation: Soviet Math. Dokl. {\bf 33}
(1986), no. 1, 227--230.

\bibitem[W]{W}
B.L. van der Waerden. {\it Ein einfaches Beispiel einer
nichtdifferenzierbaren stetigen Funktion,} Math. Zeitschr. {\bf
32} (1930), 474--475.

\bibitem[Z]{Z}
L. Zielinski. {\it Asymptotic distribution of eigenvalues for
elliptic boundary value problems,} Asymptot. Anal. {\bf 16}
(1998), no. 3-4, 181--201.

\end{thebibliography}
\end{document}